\documentclass[Journal,letterpaper,NoLineNumbers]{ascelike-new}
\usepackage[utf8]{inputenc}
\usepackage[T1]{fontenc}
\usepackage{lmodern}
\usepackage{graphicx}
\usepackage[style=base,figurename=Fig.,labelfont=bf,labelsep=period]{caption}
\usepackage{multirow}
\usepackage{subcaption}
\usepackage{amsmath}
\usepackage{newtxtext,newtxmath}
\usepackage{natbib}
\usepackage{booktabs}

\usepackage[colorlinks=true,citecolor=red,linkcolor=black]{hyperref}

\newcounter{pdfbookmarksec}

\newcommand{\bmsection}[1]{%
    \stepcounter{pdfbookmarksec}%
    \pdfbookmark[1]{#1}{bmsec.\arabic{pdfbookmarksec}}%
    \section{#1}%
}

\newcommand{\bmsectionstar}[1]{%
    \stepcounter{pdfbookmarksec}%
    \pdfbookmark[1]{#1}{bmsec.\arabic{pdfbookmarksec}}%
    \section*{#1}%
}
\NameTag{Ko, \today}
\begin{document}

\title{Quadratic Point Estimate Method for Uncertainty Quantification with Dependent Non-Gaussian Inputs}

\author[1]{Minhyeok Ko}
\author[2]{Konstantinos G. Papakonstantinou}

\affil[1]{Department of Civil and Construction Engineering and Management, The University of Texas at Tyler, Email: mko@uttyler.edu}
\affil[2]{Department of Civil and Environmental Engineering, The Pennsylvania State University}

\maketitle

\begin{abstract}
As an extension of the Point Estimate Method (PEM) to evaluate probabilistic moments of quantities of interest (QoI) in general $n$-dimensional spaces, the Quadratic Point Estimate Method (QPEM) has been recently developed. This new method is defined to fully represent up to fifth-order input moments in the Gaussian space, providing general analytical expressions for sample locations and weights, without requiring any numerical optimization. The QPEM can significantly improve the estimation accuracy of the output QoI moments, in relation to PEM-based methods whose numbers of sigma points grow linearly with the problem dimension, while at the same time having an affordable and competitive computational cost up to a considerable number of dimensions. The QPEM is further enhanced in this work by enabling copula integration into the framework, which enables effective modeling of the joint input probability density function by estimating marginals and the dependence structure of the involved  random variables. The validity and efficient performance of the copula-based QPEM are showcased against numerous other sampling methods in various examples considering two practical scenarios: (i) when the joint dependence structure can be inferred from data, and (ii) when only marginal distributions and correlation matrices are known.
\end{abstract}

\newpage
\bmsection{Introduction}
\label{sec:Introduction}
In engineering applications, uncertainty quantification (UQ) through nonlinear models frequently requires computing probabilistic moment integrals. Among available numerical methods, the Monte Carlo (MC) method is extensively employed due to its versatility. However, standard MC approaches typically demand significant computational resources when applied to computationally intensive models. To address this limitation, variance reduction techniques such as quasi-Monte Carlo (QMC) methods (\citealt{caflisch1998monte,sobol1998quasi,jank2005quasi}) and Latin Hypercube Sampling (LHS) (\citealt{mckay1979lhs,olsson2002latin,shields2016generalization}) have been developed as computationally efficient alternatives. These techniques exhibit faster convergence rates than traditional MC methods under specific conditions and dimensional constraints. Quadrature techniques have also been employed to estimate moment integrals, demonstrating efficiency and accuracy in low-dimensional cases (\citealt{PAPAKONSTANTINOU2013286}). However, these techniques are generally inefficient for probabilistic problems involving high-dimensional spaces.

The Point Estimate Method (PEM), first introduced by \cite{rosenblueth1975point}, represents a dedicated probabilistic approach with a long development history in engineering. Since its introduction, numerous PEM variations have been developed as powerful yet simple approximation methods for estimating the first few moments of output probability density function (PDF) (\citealt{harr1989,hong1998,julier1997new,julier2002scaled,julier2002spherical,tenne2003higher,xiao2018reliability,adurthi2018conjugate,papakonstantinou2022scaled,amir2022scaled}). PEM-based methods often provide accurate mean and standard deviation estimates with notable computational efficiency. However, most existing PEMs are inadequate for estimating higher-order moments such as skewness and kurtosis without substantial computational expense. The accuracy of higher-order moment estimation depends on adequately capturing corresponding high-order input moments, which typically requires an increased number of sampling points, often also called sigma points. Consequently, the main computational cost of PEMs relates to the number of sigma points utilized, creating an inherent trade-off between computational efficiency and accuracy. Additionally, many PEM methods rely on optimization procedures to determine sigma point locations and weights (\citealt{rosenblueth1975point,adurthi2018conjugate}), which complicates their practical application.

To address these limitations, the Quadratic Point Estimate Method (QPEM) has been recently introduced by \cite{ko2025quadratic}. The QPEM is explicitly designed to enhance higher-order moment estimation accuracy without significantly increasing computational demands. The method completely captures up to fifth-order input moments of Gaussian distributions without requiring computationally intensive optimization procedures for sigma point determination. As a result, the QPEM effectively reduces errors in higher-order output moment estimation, offering notable improvements over traditional PEM methods.

The accurate representation of statistical dependencies among input variables remains a fundamental challenge in UQ. Conventional frameworks frequently assume multivariate normality, primarily due to its analytical tractability and ease of implementation within existing probabilistic analysis tools. However, such Gaussian assumptions do not always hold in realistic engineering systems, where dependencies often arise from nonlinear physical interactions or heterogeneous material properties. Recent studies (\citealt{fang2013modified,radu2019reliability}) have demonstrated that these simplifying assumptions can yield biased or misleading results, especially in estimating higher-order moments critical to risk and reliability assessments. The limitations become particularly evident when analyzing systems governed by extreme or tail-dependent behavior, where non-Gaussian features dominate the probabilistic response. These observations highlight the necessity of more general frameworks capable of representing a broad spectrum of dependency structures beyond the Gaussian domain.

In this context, recent developments have sought to extend PEM frameworks to handle dependent and non-Gaussian joint inputs. A data-centric Gaussian-mixture-based PEM has been proposed to improve joint moment prediction for dependent parameters while avoiding the inaccuracies of single-Gaussian mappings (\citealt{xie2019robust}). PEM can also accommodate dependent non-Gaussian inputs through distribution-specific transformations, such as the inverse Rosenblatt mapping (\citealt{schenkendorf2014general}). Furthermore, complementary sigma-point formulations, such as the randomized Unscented Transform (UT), enhance integration accuracy through stochastic generation of sigma points, achieving asymptotically exact moment estimation for nonlinear, non-Gaussian systems while maintaining tractable computational cost (\citealt{straka2012randomized}). These studies collectively demonstrate the potential of combining advanced dependency modeling with efficient moment-based formulations, motivating the development of a unified framework that systematically integrates such capabilities.

A number of advanced methods have been developed for input characterization, including parametric and nonparametric marginal modeling (\citealt{silverman2018density,soize2000nonparametric}), mixture-based representations (\citealt{mclachlan2019finite,bishop2006pattern}), copula frameworks (\citealt{zhang2020copula,ding2023copula}), kernel-based density (\citealt{jung2024confidence}), and modern high-dimensional density models such as normalizing flows and variational autoencoders (\citealt{rezende2015variational}). Among these, copula-based approaches have gained prominence due to their ability to decouple marginal distributions from dependency structures, allowing flexible and interpretable formulations of multivariate models (\citealt{he2021dependence,torre2019general}). In particular, vine copulas provide a hierarchical and modular construction that decomposes high-dimensional dependencies into a set of bivariate copulas, thereby capturing asymmetric and nonlinear relationships with notable computational efficiency. However, integrating vine copulas with efficient moment estimation methods has remained largely unexplored, representing a significant opportunity for methodological advancement (\citealt{torre2019general,xu2020vine}). From a practical standpoint, input characterization can be broadly classified into two representative cases: (i) when observational data are available to explicitly model and calibrate the dependence structure; and (ii) when only marginal distributions and pairwise correlations are known, requiring approximate dependency representations (\citealt{kumar2023bayesian}). Distinguishing between these two scenarios is crucial, as it delineates the attainable fidelity of the probabilistic model and directly influences the accuracy of moment estimation and probabilistic predictions in the subsequent analysis. These distinctions motivate the development of a unified framework that effectively integrates flexible dependency modeling with efficient probabilistic moment computation.

This work addresses this critical gap by developing a comprehensive framework that combines the computational efficiency of QPEM with the modeling flexibility of vine copulas. We enhance and generalize the QPEM by integrating vine copulas into the framework, enabling effective modeling of joint input PDFs. Our approach constructs vine copula models that capture complex, non-Gaussian dependencies. The method maintains the computational efficiency that makes QPEM practical and desirable, and provides a general, extensible foundation for uncertainty quantification in complex engineering applications.


The remainder of this paper is organized as follows. The QPEM formulation and its deterministic sigma-point construction for probabilistic moment estimation are first reviewed. The fundamentals of copulas and vine copulas are then introduced, including Gaussian copulas, the Nataf transformation, and pair-copula constructions. Next, the proposed copula-based QPEM framework is presented under two representative levels of input information: (i) when sufficient data are available to infer the dependence structure and (ii) when only the marginal distributions and correlation matrix are known. The framework is subsequently evaluated through six numerical studies involving rock slope stability, tunnel excavation, and planar and spatial truss systems, with comparisons against Monte Carlo simulation, Latin hypercube sampling, quasi-Monte Carlo sampling, Smolyak Gauss-Hermite sparse-grid quadrature, and Hong’s point estimate method. Finally, the main findings and conclusions are summarized.

\bmsection{Quadratic Point Estimate Method}
\label{sec:QPEM}
To address effectively the propagation of uncertainty through nonlinear functions in probabilistic mechanics, the Point Estimate Method (PEM) constitutes one of the earliest and most widely adopted deterministic sampling techniques. Originally proposed by \cite{rosenblueth1975point}, the PEM approximates statistical moments of a system response by evaluating the governing model at a limited number of carefully chosen sigma points, rather than performing large-scale random sampling. In this framework, each sigma point represents a specific combination of input variable values, selected such that the corresponding weights satisfy a prescribed set of Moment Constraint Equations. By transforming multidimensional probabilistic integrations into finite weighted summations, PEMs dramatically reduce the number of required model evaluations while maintaining acceptable accuracy for the lower-order moments of interest. However, the original Rosenblueth scheme, employing a two-point expansion for each random variable, results in an exponential growth of model evaluations ($2^n$ for an $n$-dimensional problem), which inevitably becomes computationally prohibitive for moderately high dimensions.

While alternative methods, including Harr’s (\citealt{harr1989}), Hong’s (\citealt{hong1998}), and the Unscented Transformation variants (\citealt{julier1997new,julier2002scaled,julier2002spherical}), were subsequently developed to mitigate this curse of dimensionality and improve moment-matching capabilities, they still exhibit distinct shortcomings. These techniques, typically employing symmetric arrangements of sigma points with linear computational complexity (often $2n+1$ evaluations), successfully recover the first two or three input moments. Nevertheless, the confinement of their sigma points to the principal axes of the standardized variable space prevents the accurate estimation of cross-moments among variables. This limitation is critical, as these quantities become particularly influential in nonlinear, coupled systems. Consequently, existing PEMs face an inherent trade-off between computational cost and accuracy, especially when higher-order or mixed probabilistic moments are sought.

\begin{figure}[t]
    \centering
    \includegraphics[width=0.5\textwidth]{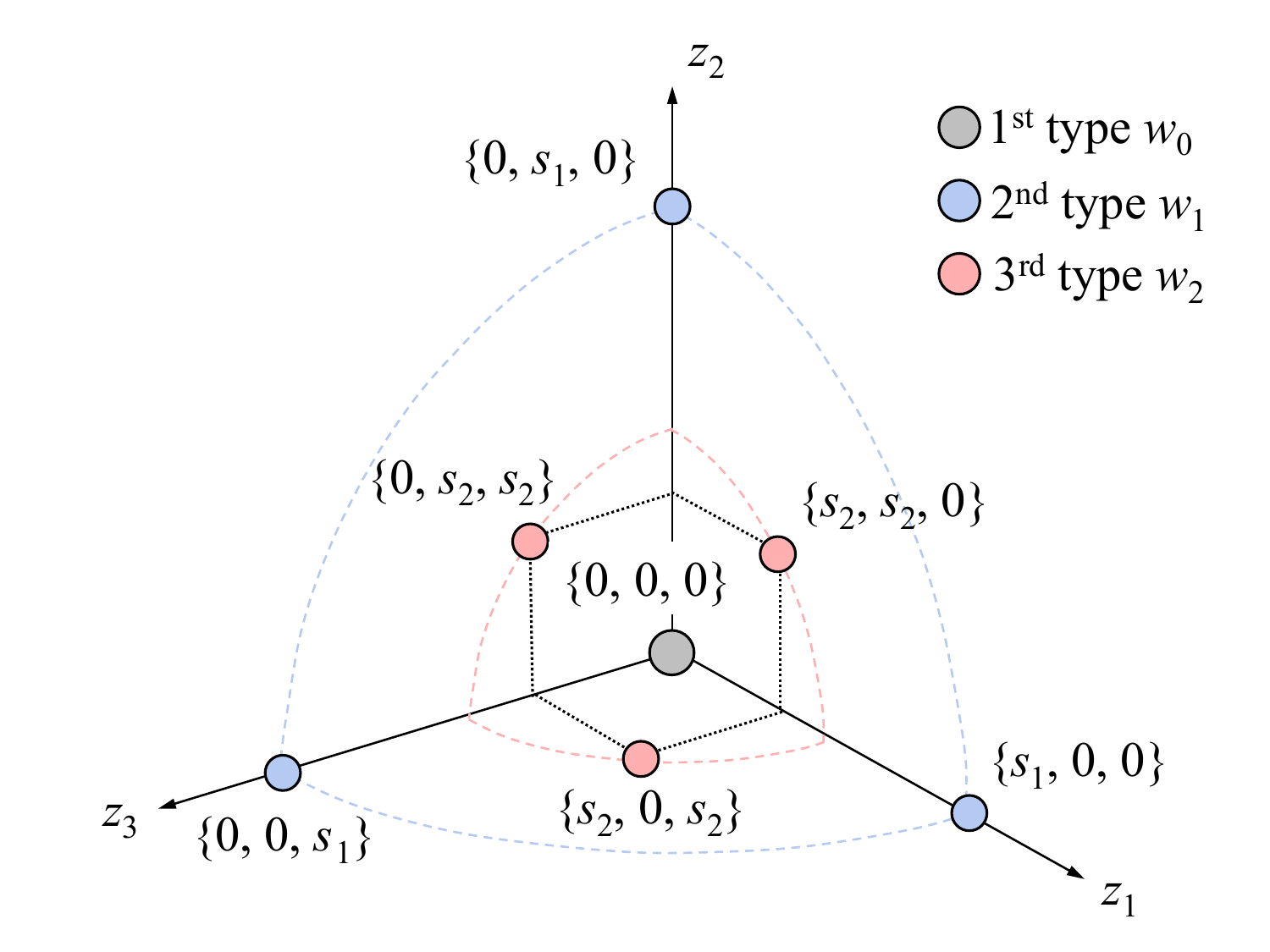}
\caption{Schematic representation of QPEM points in 3-$\mathbb{D}$}
\label{fig:qpem}
\end{figure}

To overcome these critical limitations, the Quadratic Point Estimate Method (QPEM) (\citealt{ko2025quadratic}) has recently been introduced. The QPEM systematically extends PEM accuracy to higher-order statistics without relying on optimization procedures or exponentially growing sampling sets. The QPEM framework is distinguished by three key attributes: (a) fully symmetric sigma points ensuring automatic satisfaction of odd-order zero moments; (b) analytical, closed-form expressions for the sigma point coordinates and associated weights, effectively eliminating the need for cumbersome optimization procedures; and (c) enhanced accuracy in estimating higher-order output moments based on $2n^2+1$ sigma points evaluations. Specifically, the QPEM is capable of capturing exactly up to the fifth-order moments of an $n$-dimensional standard Gaussian random variable $\mathbf{z}$. Figure~\ref{fig:qpem} schematically illustrates the fully symmetric arrangement of the sigma points, classified according to their spatial configuration and associated weights.

The sigma point set is comprised of: (a) a single point at the origin $\mathbf{0}$ with weight $w_0$; (b) $2n$ points positioned along each orthogonal axis at a distance $s_1$ from the origin, each with weight $w_1$; and (c) $2n(n-1)$ points generated through all permutations and sign variations of the coordinates $\{s_2, s_2, 0, \dots, 0\}^T$, each with weight $w_2$. The $2n^2 + 1$ employed sigma points thus enable the complete capture of fifth-order input moments. Detailed theoretical derivations of these sigma points and corresponding weights are provided in the original QPEM formulation (\citealt{ko2025quadratic}).

The analytical expressions defining the sigma point locations and associated weights are given by:

\begin{equation}
	\label{eq:qpem}
	\left\{\begin{array}{ll}
            {{s}_{0}} &=0\\
		{{s}_{1}} &=r\\ 
		{{s}_{2}} &=\left[\dfrac{r^2\left(n-1\right)}{r^2+n-4}\right]^{1/2}
	\end{array} \right. %
	\left\{ \begin{array}{ll}
            {{w}_{0}} &= 1-2nw_1-2n\left(n-1\right)w_2 \\
		{{w}_{1}} &= \dfrac{4-n}{2r^4} \\
		{{w}_{2}} &= \dfrac{1}{4}\left[\dfrac{r^2+n-4}{r^2\left(n-1\right)}\right]^2
	\end{array} \right.
\end{equation}
where $r$ is a user-defined tuning parameter constrained by $r > \sqrt{2}$ to ensure numerical stability. The selection of $r$ influences the residual error in moment estimations beyond the fifth order. Following the recommendations in the original QPEM formulation (\citealt{ko2025quadratic}), a value of $r = 3$ is utilized in this study, as it provides a favorable balance between numerical stability and overall estimation accuracy.

Furthermore, the QPEM framework offers the flexibility to capture input moments beyond the fifth order by adjusting the scaling parameters $\zeta$ and $\xi$. These parameters effectively modify the weight assigned to the origin, thereby fine-tuning output skewness and kurtosis estimation without incurring additional computational costs. In accordance with the default recommendations provided in the original QPEM formulation (\citealt{ko2025quadratic}), the values of $\zeta = -8$ and $\xi = 60$ are adopted in this study to ensure robust performance and consistency with the original formulation. Table~\ref{table:qpem_points_weights} summarizes the sigma points and weights for the QPEM formulation used in this work.

Irrespective of the original distribution, sigma points ($\mathbf{S}_i$) are initially identified in the standard normal space ($z$-space) and subsequently transformed to the original variable space ($x$-space) via an appropriate transformation $\mathbf{X}_i=\mathcal{T}(\mathbf{S}_i)$. The response at each transformed sigma point is evaluated through the computational model $\mathbf{Y}_i=\mathbf{M}(\mathbf{X}_i)$. The first four probabilistic response moments are computed as follows:

\begin{align}
    \label{eq:qpem_mean_equation}
    E\left[ \mathbf{y} \right]\equiv \mathbf{\bar{y}} & =\displaystyle \sum\limits_{i=0}^{2{{n}^{2}}}{W_{i}^{\left( 1 \right)}{{\mathbf{Y}}_{i}}} \\ 
    \label{eq:qpem_moments_equation}
    E\left[ {{\left( \mathbf{y}-\mathbf{\bar{y}} \right)}^{k}} \right] & =\displaystyle \sum\limits_{i=0}^{2{{n}^{2}}}{W_{i}^{\left( k \right)}\left(\mathbf{Y}_i-\overline{\mathbf{y}}\right)^{k}}
\end{align}
where $k=2,\,3,\,4$ denotes the order of the central moment.

\begin{table*}[ht!]
	\caption{Sigma points and weights selection for QPEM (\citealp{ko2025quadratic})}\label{table:qpem_points_weights}
	\begin{tabular*}{\textwidth}[t]{@{}lll@{}}
		\hline
		\begin{tabular}[t]{l} \textbf{Weights}\end{tabular}&  & \begin{tabular}[t]{|l} 
			$W_{0}^{(1)}=W_{0}^{(2)}=w_0$   \quad $W_{0}^{(3)}=w_0+\zeta$  \quad $W_{0}^{(4)}=w_0+\xi$\\
			$W_i^{(k)\,\,a} = \begin{cases}
				w_1, & \text{for $i=1,\cdots, 2n$} \\ 
				w_2, & \text{for $i=2n+1,\cdots, 2n^2$} 
			\end{cases}$ \end{tabular} \\ 
		\hline
		\multirow{2}{*}{\begin{tabular}[t]{l} \textbf{Sigma}\\ \textbf{Point} \end{tabular}} & 1st type & \begin{tabular}[t]{|l} 
			$\mathbf{S}_0=\mathbf{0}_{n \times 1}$
		\end{tabular}\\
		\begin{tabular}[t]{l}	\end{tabular}							& 2nd type & \begin{tabular}[t]{|lll}
			$\mathbf{S}_1=\{s_1, 0,\cdots, 0\}^T$& $\cdots$ & $\mathbf{S}_n=\{0,\cdots, 0, s_1\}^T$\\
			$\mathbf{S}_{n+1}=\{-s_1, 0,\cdots, 0\}^T$& $\cdots$ & $\mathbf{S}_{2n}=\{0,\cdots, 0, -s_1\}^T$\end{tabular} \\
		\begin{tabular}[t]{l}	\end{tabular}							& 3rd type$^{b}$ &  \begin{tabular}[t]{|lll}
			$\mathbf{S}_{2n+1}=\{s_2, s_2, 0, \cdots, 0\}^T$ & $\cdots$ & $\mathbf{S}_{2n^2-3}=\{0,\cdots, 0, s_2, s_2\}^T$\\
			$\mathbf{S}_{2n+2}=\{-s_2, s_2, 0,\cdots, 0\}^T$& $\cdots$ & $\mathbf{S}_{2n^2-2}=\{0,\cdots, 0, -s_2, s_2\}^T$\\
			$\mathbf{S}_{2n+3}=\{s_2, -s_2, 0,\cdots, 0\}^T$& $\cdots$ & $\mathbf{S}_{2n^2-1}=\{0,\cdots, 0, s_2, -s_2\}^T$\\
			$\mathbf{S}_{2n+4}=\{-s_2, -s_2, 0,\cdots, 0\}^T$& $\cdots$ & $\mathbf{S}_{2n^2}=\{0,\cdots, 0, -s_2, -s_2\}^T$\end{tabular} \\ 
		\hline
		\multicolumn{3}{l}{\begin{tabular}{@{}l@{}} $^a$ Associated weights for calculating $k^{th}$ moments of the response function\\$^b$ Obtained by the permutation and the change of sign of $\mathbf{S}_{2n+1}$ \end{tabular}}\\ 
	\end{tabular*}
\end{table*}

\bmsection{Copulas and vine copulas}
\label{sec:copula}
\subsection{Copula and Sklar's theorem}
\begin{equation}
    \label{eq:copula_decom1}
    F\left(\mathbf{x}\right) = C_\mathbf{x}\left(F_1(x_1), \cdots, F_n(x_n)\right)
\end{equation}
where $F$ denotes the joint CDF, and $F_1, \dots, F_n$ represent the corresponding univariate marginal CDFs. By differentiating Eq.~\eqref{eq:copula_decom1}, the joint probability density function (PDF) $f(\mathbf{x})$ is derived via the chain rule as:
\begin{equation}
\label{eq:pdf_chain_rule}
\begin{split}
f(\mathbf{x})&=\dfrac{\partial^n F(\mathbf{x})}{\partial x_1 \ldots \partial x_n}\\
&=c_{\mathbf{x}}\left\{F_1\left(x_1\right), \ldots, F_n\left(x_n\right)\right\} \cdot \displaystyle  \prod_{i=1}^n f_i\left(x_i\right)
\end{split}
\end{equation}
where $c_{\mathbf{x}}(\mathbf{u}) = {\partial^n C_{\mathbf{x}}(\mathbf{u})}/{\partial u_1 \dots \partial u_n}$ is the copula density function, and $f_i(x_i)$ are the marginal PDFs. This formulation highlights the advantage of copulas: the separation of the dependence structure ($c_{\mathbf{x}}$) from the marginal behaviors ($f_i$).

In the context of copula-based dependency modeling, rank correlation coefficients are widely adopted for quantifying correlation, addressing shortcomings of Pearson's linear correlation coefficient ($\rho$) in characterizing nonlinear dependencies. Specifically, Kendall's rank correlation coefficient ($\tau$) is particularly favored in this work, as it possesses superior attributes in capturing nonlinear relationships independent of the marginal distributions. Furthermore, unlike linear correlation, Kendall’s $\tau$ remains invariant under strictly monotonic transformations of marginals and is explicitly linked to the copula function by \cite{nelsen2007copula}:
\begin{equation}
    \tau=4 \int_{0}^1 \int_{0}^1 C\left(u_1, u_2\right) d C\left(u_1, u_2\right)-1
\end{equation}

\subsection{Elliptical copulas and Nataf transformation}
Within the broad spectrum of copula families, elliptical copulas such as the Gaussian and Student-$t$ are predominantly utilized in engineering applications largely due to their mathematical tractability and ease of simulation. The Gaussian copula, $C_{G}$, is characterized by the multivariate standard normal distribution $\Phi_n$ and a correlation matrix $\mathbf{R}$. For a vector $\mathbf{u} \in [0,1]^n$, the Gaussian copula is formulated as:
\begin{equation}
    \label{eq:gaussian_copula}
    C_{G}(\mathbf{u} \mid \mathbf{R}) = \Phi_n\left(\Phi^{-1}(u_1), \dots, \Phi^{-1}(u_n) \mid \mathbf{R}\right)
\end{equation}
where $\Phi^{-1}$ denotes the inverse CDF of the standard normal distribution.

This formulation is intrinsically linked to the Nataf transformation, a technique extensively applied in structural reliability frameworks (\citealt{der1986structural,lebrun2009generalization}). The Nataf model operates by transforming a correlated random vector $\mathbf{X}$ into a standard normal vector $\mathbf{Z}$ via the marginal transformation $z_i = \Phi^{-1}(F_i(x_i))$. Consequently, the joint PDF of $\mathbf{X}$ within the Nataf framework is derived as:
\begin{equation}
    \label{eq:nataf_pdf}
    f(\mathbf{x}) = \phi_n(\mathbf{z}, \mathbf{R}') \prod_{i=1}^n \frac{f_i(x_i)}{\phi(z_i)}
\end{equation}
where $\phi_n(\cdot, \mathbf{R}')$ represents the $n$-dimensional standard normal density characterized by the correlation matrix $\mathbf{R}'$, and $\phi(\cdot)$ is the univariate standard normal density.

A fundamental attribute of the Nataf transformation lies in its handling of correlation. The correlation coefficient $\rho'_{ij}$ in the Gaussian space (elements of $\mathbf{R}'$) is related to the correlation coefficient $\rho_{ij}$ in the physical space through the integral equation:
\begin{equation}
    \label{eq:nataf_correlation}
    \rho_{ij} = \int_{-\infty}^{\infty} \int_{-\infty}^{\infty} \left(\frac{x_i-\mu_i}{\sigma_i}\right) \left(\frac{x_j-\mu_j}{\sigma_j}\right) \phi_2(z_i, z_j, \rho'_{ij}) \, dz_i \, dz_j
\end{equation}
This relationship effectively demonstrates the mathematical equivalence between the Nataf transformation and the modeling of the dependency structure using a Gaussian copula (\citealt{lebrun2009generalization}).

Despite its widespread adoption, the Gaussian copula (and by extension the Nataf transformation) exhibits distinct shortcomings regarding the flexibility of the dependence structure. Since it is exclusively defined by the linear correlation matrix $\mathbf{R}$, it is inherently restricted to capturing only linear relationships between variables. Real-world stochastic phenomena may often exhibit complex, non-linear dependencies that a simple linear correlation coefficient may fail to represent adequately. Consequently, the Gaussian model can prove essentially too rigid for high-dimensional applications in such settings where accurate modeling of these complex interactions is imperative, thus motivating the adoption of vine copulas as discussed in the following section.

\subsection{Vine copulas}
As the input dimension $n$ increases, the accurate representation of dependencies among variables via an $n$-dimensional multivariate copula presents significant challenges. Although multivariate copula families exist, they frequently fail to adequately capture complex dependencies (\citealt{aas2009pair}). To circumvent these limitations, \cite{bedford2002vine} introduced a flexible methodology known as Pair-Copula Construction (PCC), which systematically constructs any $n$-dimensional copula density from products of bivariate copulas (pair-copulas).

This framework leverages the fundamental definition of conditional densities, given by:
\begin{equation}
\label{eq:def_conditional_density}
f_{i|1,\ldots,i-1}(x_i \mid x_1, \dots, x_{i-1})=\frac{f(x_1, \dots, x_{i-1}, x_i)}{f(x_1, \dots, x_{i-1})}
\end{equation}
Based on Eq.~\eqref{eq:def_conditional_density}, the joint PDF $f(\mathbf{x})$ can be factorized as:
\begin{equation}
\label{eq:pdf_factorization}
f(\mathbf{x})=f_{1}(x_1) f_{2|1}(x_2 \mid x_1) \cdots f_{n|1,2,\ldots,n-1}(x_n \mid x_1, \dots, x_{n-1})
\end{equation}
Furthermore, each term in Eq.~\eqref{eq:pdf_factorization} admits further decomposition utilizing pair-copulas and conditional marginal densities:
\begin{equation}
\label{eq:condi_pdf_copula}
f_{i \mid \mathbf{i}}(x_i \mid \mathbf{v})=c_{ij \mid \mathbf{i}_{\overline{j}}}\left\{F_{i \mid \mathbf{i}_{\overline{j}}}\left(x_i \mid \mathbf{v}_{\overline{j}}\right), F_{j \mid \mathbf{i}_{\overline{j}}}\left(x_j \mid \mathbf{v}_{\overline{j}}\right)\right\}f_{i\mid \mathbf{i}_{\overline{j}}}\left(x_i \mid \mathbf{v}_{\overline{j}}\right)
\end{equation}
where $x_j$ is arbitrarily excluded from the vector $\mathbf{v}$, and $\mathbf{v}_{\overline{j}}$ denotes the remaining elements. The indices $\mathbf{i}$ and $\mathbf{i}_{\overline{j}}$ correspond to the vectors $\mathbf{v}$ and $\mathbf{v}_{\overline{j}}$, respectively. By iteratively combining Eqs.~\eqref{eq:pdf_factorization} and \eqref{eq:condi_pdf_copula}, a factorization of $f(\mathbf{x})$ can be achieved employing exclusively marginal distributions and pair-copulas.

\begin{figure}[!t]
    \centering
    \includegraphics[width=0.6\textwidth]{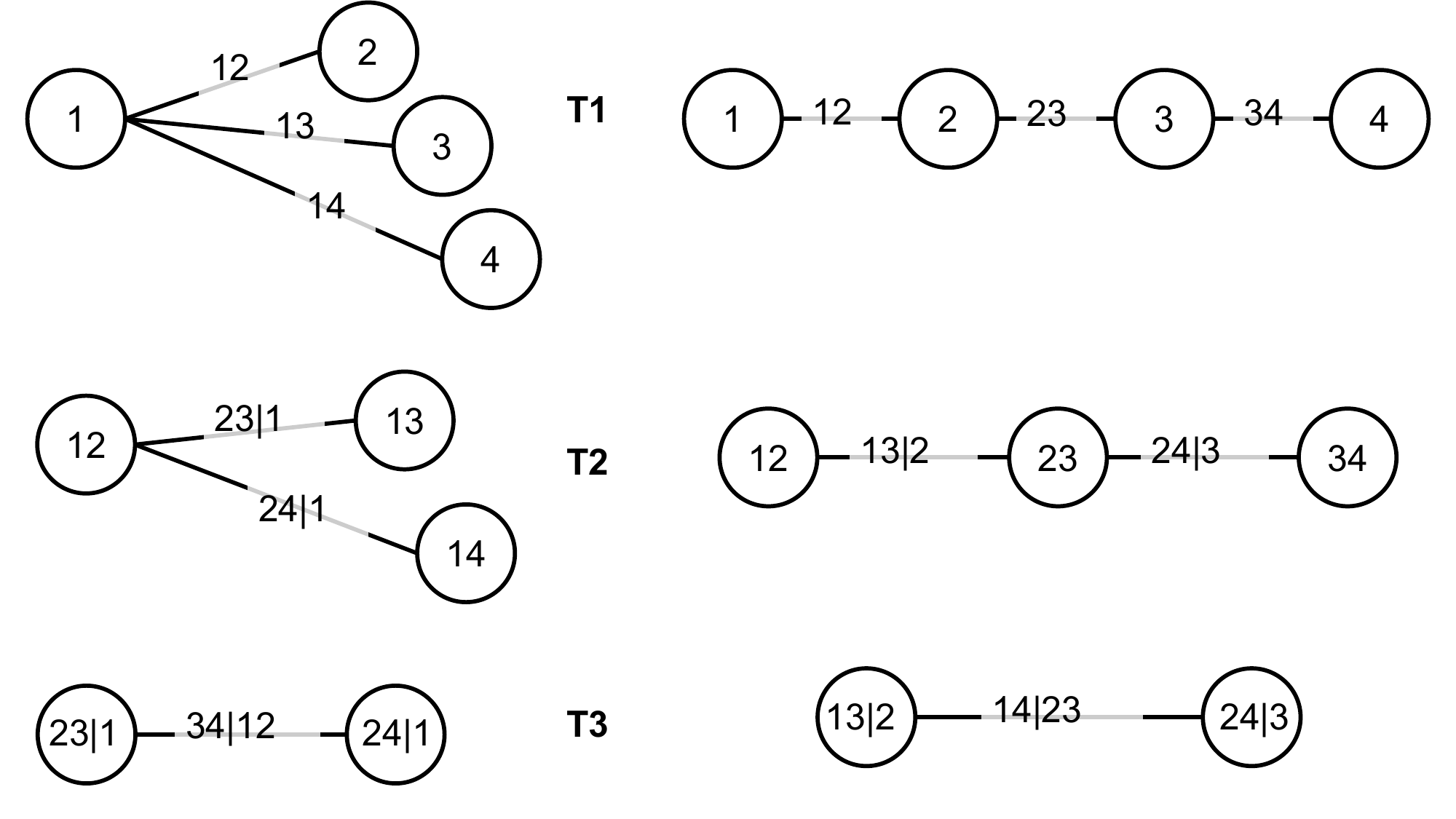}
    \caption{Representative C-vine (left) and D-vine (right) in a four-dimensional space.}
\label{fig:c_d_vine}
\end{figure}

To illustrate this concept, consider a three-dimensional case ($n=3$). The conditional PDF $f_{3 \mid 12}\left(x_3 \mid x_1, x_2\right)$ can be expressed as:
\begin{equation}
    \label{eq:3d_cond_pdf}
    \begin{split}
        f_{3 \mid 12}\left(x_3 \mid x_1, x_2\right) & =\frac{f_{13 \mid 2}\left(x_1, x_3 \mid x_2\right)}{f_{1 \mid 2}\left(x_1 \mid x_2\right)}\\
        &=\frac{c_{13 \mid 2}\left(F_{1 \mid 2}\left(x_1 \mid x_2\right), F_{3 \mid 2}\left(x_3 \mid x_2\right)\right) \cdot f_{1 \mid 2}\left(x_1 \mid x_2\right) \cdot f_{3 \mid 2}\left(x_3 \mid x_2\right)}{f_{1 \mid 2}\left(x_1 \mid x_2\right)} \\
        & =c_{13 \mid 2}\left(F_{1 \mid 2}\left(x_1 \mid x_2\right), F_{3 \mid 2}\left(x_3 \mid x_2\right)\right) \cdot f_{3 \mid 2}\left(x_3 \mid x_2\right) \\
        & =c_{13 \mid 2}\left(F_{1 \mid 2}\left(x_1 \mid x_2\right), F_{3 \mid 2}\left(x_3 \mid x_2\right)\right) \cdot c_{23}\left(F_2\left(x_2\right), F_3\left(x_3\right)\right) \cdot f_3\left(x_3\right)
    \end{split}
\end{equation}
with the conditional CDFs given by:
\begin{equation}
    \label{eq:3d_cond_cdf}
    \begin{split}
    & F_{1 \mid 2}\left(x_1 \mid x_2\right)=C_{1 \mid 2}\left\{F_1\left(x_1\right) \mid F_2\left(x_2\right)\right\}=\frac{\partial C_{12}\left\{F_1\left(x_1\right), F_2\left(x_2\right)\right\}}{\partial F_2\left(x_2\right)} \\
    & F_{3 \mid 2}\left(x_3 \mid x_2\right)=C_{3 \mid 2}\left\{F_3\left(x_3\right) \mid F_2\left(x_2\right)\right\}=\frac{\partial C_{32}\left\{F_3\left(x_3\right), F_2\left(x_2\right)\right\}}{\partial F_2\left(x_2\right)}
    \end{split}
\end{equation}
Consequently, from Eqs.~\eqref{eq:3d_cond_pdf} and \eqref{eq:3d_cond_cdf}, the joint PDF $f\left(x_1, x_2, x_3\right)$ can be constructed as:
\begin{equation}
    \begin{split}
    f\left(x_1, x_2, x_3\right) & =f_1\left(x_1\right) \cdot f_{2 \mid 1}\left(x_2 \mid x_1\right) \cdot f_{3 \mid 12}\left(x_3 \mid x_1, x_2\right) \\
    & =f_1\left(x_1\right) \cdot f_2\left(x_2\right) \cdot f_3\left(x_3\right) \cdot c_{12}\left\{F_1\left(x_1\right), F_2\left(x_2\right)\right\} \\
    & \quad \cdot c_{23}\left\{F_2\left(x_2\right), F_3\left(x_3\right)\right\} \cdot c_{13 \mid 2}\left\{F_{1 \mid 2}\left(x_1 \mid x_2\right), F_{3 \mid 2}\left(x_3 \mid x_2\right)\right\}
    \end{split}
\end{equation}
Various decompositions of $f(\mathbf{x})$ are attainable due to the dependence of pair-copula selection on the conditioning variables at each step. For high-dimensional problems, the number of potential PCCs expands combinatorially. To manage this complexity, \cite{bedford2002vine} introduced the concept of a \textit{regular vine} (so-called R-vine), a graphical modeling framework (Figure~\ref{fig:c_d_vine}) that structures dependencies through a sequence of nested trees $\mathcal{V}=\left(T_1, \ldots, T_{n-1}\right)$, satisfying the following conditions (\citealt{joe2014dependence}):
\begin{itemize}
    \item[(a)] $T_1$ is a tree with nodes $N_1=\{1, \ldots, n\}$ and edge set $E_1$.
    \item[(b)] For $i=2, \ldots, n-1, T_i$ is a tree with nodes $N_i=E_{i-1}$ and edge set $E_i$.
    \item[(c)] (Proximity condition) If $\{\mathbf{a}, \mathbf{b}\} \in E_i$ where $\mathbf{a}=\left\{a_1, a_2\right\}$ and $\mathbf{b}=\left\{b_1, b_2\right\}$, then \#($\mathbf{a} \cap \mathbf{b})=1$, where \# denotes the cardinality of a set.
\end{itemize}

In this framework, the edges in tree $T_{j}$ become nodes in tree $T_{j+1}$, and the proximity condition ensures consistent conditioning. R-vines are extensively adopted to construct joint distributions via conditional pair-copulas arranged in a tree structure, yielding a class of models known as \textit{vine copulas}.
\par

Two predominant classes of R-vines are the drawable vine (D-vine) and the canonical vine (C-vine), as representatively shown in Figure~\ref{fig:c_d_vine}. For instance, the four-dimensional joint PDF corresponding to the C-vine in Figure~\ref{fig:c_d_vine} is given by:
\begin{equation}
\label{eq:4d_c_vine_decomp}
\begin{split}
    f(\mathbf{x})= & f_1\left(x_1\right)\cdot f_2\left(x_2\right)\cdot f_3\left(x_3\right)\cdot f_4\left(x_4\right) \\
    &\cdot  c_{1 2}\left\{F_1\left(x_1\right),F_2\left(x_2\right)\right\}\cdot c_{1 3}\left\{F_1\left(x_1\right),F_3\left(x_3\right)\right\} \cdot c_{1 4}\left\{F_1\left(x_1\right),F_4\left(x_4\right)\right\} \\
    &\cdot c_{2 3 \mid 1}\left\{F_{2 \mid 1}\left(x_2 \mid x_1\right), F_{3 \mid 1}\left(x_3 \mid x_1\right)\right\}\cdot c_{2 4 \mid 1}\left\{F_{2 \mid 1}\left(x_2 \mid x_1\right), F_{4 \mid 1}\left(x_4 \mid x_1\right)\right\}\\
    &\cdot c_{3 4 \mid 1 2}\left\{F_{3 \mid 12}\left(x_3 \mid x_1, x_2\right), F_{4 \mid 12}\left(x_4 \mid x_1, x_2\right)\right\}
\end{split}
\end{equation}
while the joint PDF corresponding to the D-vine in Figure~\ref{fig:c_d_vine} is expressed as:
\begin{equation}
\label{eq:4d_d_vine_decomp}
\begin{split}
    f(\mathbf{x})= & f_1\left(x_1\right)\cdot f_2\left(x_2\right)\cdot f_3\left(x_3\right)\cdot f_4\left(x_4\right) \\
    &\cdot  c_{1 2}\left\{F_1\left(x_1\right),F_2\left(x_2\right)\right\}\cdot c_{2 3}\left\{F_2\left(x_2\right),F_3\left(x_3\right)\right\} \cdot c_{3 4}\left\{F_3\left(x_3\right),F_4\left(x_4\right)\right\} \\
    &\cdot c_{1 3|2}\left\{F_{1 \mid 2}\left(x_1 \mid x_2\right), F_{3 \mid 2}\left(x_3 \mid x_2\right)\right\}\cdot c_{2 4|3}\left\{F_{2 \mid 3}\left(x_2 \mid x_3\right), F_{4 \mid 3}\left(x_4 \mid x_3\right)\right\}\\
    &\cdot c_{1 4|2 3}\left\{F_{1 \mid 23}\left(x_1 \mid x_2, x_3\right), F_{4 \mid 23}\left(x_4 \mid x_2, x_3\right)\right\}
\end{split}
\end{equation}

\bmsection{Copula-based QPEM}
\label{sec:copula-based QPEM}
\subsection{Rosenblatt transformation supported by vine copulas}
As discussed earlier, the sigma points $\mathbf{S}_i$ in the QPEM are initially generated in the standard normal space and then transformed into the original input $x$-space via $\mathbf{X}_i=\mathcal{T}(\mathbf{S}_i)$, where $\mathcal{T}$ denotes the mapping function. The Rosenblatt transformation is a general technique that maps a vector of dependent random variables $\mathbf{x}$ to a vector $\mathbf{z}$ with independent components. Its inverse enables the reconstruction of dependent variables from independent ones. 

A major limitation of the Rosenblatt transformation is that it requires full knowledge of the joint PDF or a sequence of conditional PDFs, which can often be unavailable in practice. To overcome this, vine copulas offer a powerful framework to construct the joint PDF, even in high-dimensional and complex dependency settings. When supported by vine copulas, the Rosenblatt transformation becomes feasible and flexible for use in QPEM, enabling accurate representation of variable dependencies.

The Rosenblatt transformation can be expressed in terms of pair-copulas using C-vine structure as follows:
\begin{equation}
    \begin{split}
        & u_1=\Phi\left(z_1\right)=F_1\left(x_1\right) \\
        & u_2=\Phi\left(z_2\right)=h\left(F_2\left(x_2\right), F_1\left(x_1\right)\right) \\
        & u_3=\Phi\left(z_3\right)=h\left(h\left(F_3\left(x_3\right), F_1\left(x_1\right)\right), h\left(F_2\left(x_2\right), F_1\left(x_1\right)\right)\right)\\
        & \cdots
    \end{split}
\end{equation}
where $\mathbf{u}$ = $\left\{u_1, \ldots,u_n\right\}$ are independent uniform random variables on $[0,1]$, and $\mathbf{z}$ = $\left\{z_1, \ldots, z_n\right\}$ are independent standard Gaussian random variables. $h(p,q)$ denotes the conditional distribution function derived from a bivariate copula $C(p,q)$, formally defined as:
\begin{equation}
    \label{eq:h-function}
    h(p,q)=\frac{\partial C(p,q)}{\partial q}
\end{equation}
This function allows one to express the conditional distribution of one variable given another in terms of their pair-copula. Its inverse $h^{-1}$ is used during the backward transformation from independent to dependent space. Therefore, the inverse Rosenblatt transformation maps independent Gaussian variables back to the dependent $x$-space and can be written as:
\begin{equation}
    \label{eq:inv_rosenblatt}
    \begin{split}
        & x_1=F_1^{-1}\left(u_1\right) \\
        & x_2=F_2^{-1}\left(h^{-1}\left(u_2, F_1\left(x_1\right)\right)\right) \\
        & x_3=F_3^{-1}\left(h^{-1}\left(h^{-1}\left(u_3, h\left(F_2\left(x_2\right), F_1\left(x_1\right)\right)\right), F_1\left(x_1\right)\right)\right)\\
        & \cdots
    \end{split}
\end{equation}
where $F_{i}^{-1}$ is the inverse function of the marginal CDF of $x_i$. A similar formulation can be derived for D-vine structures. Notably, the Nataf transformation is a special case of the Rosenblatt transformation when all pair-copulas are Gaussian (\citealt{lebrun2009generalization}).
\par
Selecting an appropriate vine structure and corresponding pair-copulas based on available data is essential, as multiple valid decompositions of the joint PDF are generally possible. Once constructed, the vine-based Rosenblatt transformation enables efficient mapping of sigma points from the standard normal space to the original input space, thereby facilitating uncertainty quantification in the presence of complex dependencies.

\subsection{Construction of input PDF using vine copula}
The construction of an appropriate vine copula model to describe the dependencies among the inputs involves three key steps: (a) selecting an appropriate R-vine structure; (b) choosing the parametric families of each pair-copula; and (c) estimating the copula parameters. Steps (a) and (b) address the representation problem by specifying a parametric model of the joint input dependencies. Step (c) finalizes the model by determining the values of the pair-copula parameters. The feasibility and accuracy of these steps are generally directly influenced by the availability and quality of input data. When sufficient data are available, it is possible to infer the dependence structure among input variables and construct a full joint PDF using vine copulas. However, in many practical settings, obtaining adequate data to fully characterize input dependencies is challenging due to resource constraints, cost, or data sparsity. To account for varying levels of available information, this study considers two distinct cases: (1) the input PDF is estimated directly from available data, and (2) the input PDF is only partially known, with information limited to marginal distributions and the covariance matrix. Each case is addressed separately in the following subsections.

\subsubsection{Case (1): Input PDF can be estimated from available data}
\label{sec:ch4-case1}
This section addresses the scenario wherein the available dataset is sufficiently informative to facilitate the direct modeling of the input distribution via vine copulas. In theory, the inference of a globally optimal vine copula model could be pursued by exhaustively evaluating the entire space of possible vine structures, copula families, and associated parameters. However, the number of potential R-vine structures exhibits an explosive growth with respect to the input dimension $n$, rendering the determination of the globally best-fitted vine copulas computationally prohibitive. Specifically, the number of distinct regular vines on $n$ nodes is given by ${n!}/2 \times {{2}^{\binom{n-2}{2}}}$ (\citealt{morales2010counting}). To illustrate the magnitude of this combinatorial conundrum, for a system with merely $n$=10 variables, the number of possible R-vine tree structures exceeds $5 \times 10^{14}$. Even if the search is constrained strictly to $\mathrm{C}$- and $\mathrm{D}$-vines, the number of possible configurations approaches approximately 2 million.

To circumvent this computational impediment, the adoption of generally efficient heuristic strategies becomes imperative. \cite{dissmann2013} proposed a sequential heuristic algorithm for the structural construction of R-vines and the selection of pair-copula families. This methodology is designed to optimize the trade-off between computational cost and estimation accuracy by prioritizing the capture of the strongest dependencies in the initial stages of the modeling process. In this framework, the strength of dependence between variable pairs is quantified using Kendall’s $\tau$, selected for its robustness and invariance to marginal distributions.

The procedure commences with the construction of the first tree $T_1$ utilizing a Maximum Spanning Tree (MST) algorithm (see Appendix~\ref{sec:app-mst} for details), where edge weights are defined by the absolute values of the pairwise empirical Kendall's $\tau$. Following the establishment of $T_1$, the optimal copula families and their corresponding parameters are identified for each edge through model selection criteria, such as the Akaike Information Criterion (AIC) or Bayesian Information Criterion (BIC). For the subsequent trees ($T_2$, $T_3$, etc.), the MST algorithm is iteratively applied to the set of edges permitted by the proximity condition, utilizing pseudo-observations transformed via the fitted copulas from the preceding trees. This sequential optimization continues until the full R-vine structure is fully determined. While this heuristic approach does not theoretically guarantee a globally optimal model (e.g., in terms of global maximum likelihood), it has been demonstrated to yield highly satisfactory results in practical engineering applications (\citealt{dissmann2013}).

Upon the rigorous specification of the R-vine model, the inverse Rosenblatt transformation, grounded in the fitted vine copulas, is utilized to map the sigma points of the QPEM, initially generated in the independent standard normal space, into the original dependent stochastic space. The comprehensive algorithm for the copula-based QPEM, applicable when adequate data are available, is summarized in Table~\ref{table:seq_method}, and the pair-copula families considered in this study are summarized in Table~\ref{table:copula_type}. For completeness, the explicit mathematical definitions of these copula families and their rotated variations are provided in Appendix~\ref{sec:app-copula-def}.

\begin{table}[!t]
	\caption{Algorithm of copula-based QPEM for adequate data scenarios}\label{table:seq_method}
	\begin{tabular}[t]{p{0.97\textwidth}}
        \hline
	\textbf{Input:} Data $\left(x_{\ell 1}, \ldots x_{\ell n}\right), \ell=1, \ldots, N$ (independent realizations from the same joint distribution).\\
        \textbf{Output:} Probabilistic output moments of nonlinear computational models.\\ [-20pt]
        \begin{itemize}
            \item[\bf{1}:] Determine regular vine copula specification.
            \item[$\qquad$\bf{1.1}:] Identify the first tree $T_1$ that maximizes the sum of absolute empirical Kendall’s $\tau$ utilizing the Maximum Spanning Tree (MST) algorithm (\citealt{cormen2022}) (see \ref{sec:app-mst}).
            \item[$\qquad$\bf{1.2}:] Select the optimal pair-copula families for edges in $T_1$ by minimizing the Akaike Information Criterion (AIC) via the Stepwise Semiparametric (SSP) estimator (\citealt{haff2013}) (see \ref{app:estimation}).
            \item[$\qquad$\bf{1.3}:] Apply MST to the graph comprising all nodes of $T_2$ with edges permissible under the proximity condition.
            \item[$\qquad$\bf{1.4}:] Proceed sequentially for the remaining trees $T_3, \dots, T_{n-1}$.
            \item[\bf{2}:] Fit and select marginal distributions based on the provided dataset.
            \item[\bf{3}:] Construct the joint PDF integrating the selected vine copulas and marginals.
            \item[\bf{4}:] Generate sigma points and weights in the standard normal space ($z$-space) via the QPEM formulation [Eq.~\eqref{eq:qpem}].
            \item[\bf{5}:] Transform sigma points from $z$-space into the original $x$-space utilizing the inverse Rosenblatt transformation defined by the vine copula model from Steps 1 and 2.
            \item[\bf{6}:] Compute the probabilistic moments of the nonlinear model response [Eqs.~\eqref{eq:qpem_mean_equation} and \eqref{eq:qpem_moments_equation}].
        \end{itemize}\\ 
        \hline
	\end{tabular}
\end{table}

\begin{table*}[t!]
\caption{Pair-copula families used in the study (\citealt{nelsen2007copula,joe2014dependence, torre2019general})}\label{table:copula_type}
\begin{tabular}{p{0.22\linewidth}p{0.19\linewidth}|p{0.3\linewidth}p{0.19\linewidth}}
\hline
Copula  & Parameter range & Copula  & Parameter range \\ \hline
Gauss   & $\theta \in (-1,1)$ & Ali-Mikhail-Haq & $\theta \in [-1,1]$  \\
t       & $\theta \in(-1,1)$, $\nu>1$ & Farlie-Gumbel-Morgenstern & $\theta \in [-1,1] \backslash\{0\}$ \\
Clayton & $\theta>0$ & Plackett & $\theta>0$ \\
Gumbel  &  $\theta \in[1,+ \inf)$ & Joe &  $\theta \in[1,+ \inf)$ \\
Frank   & $\theta \in \mathbb{R} \backslash\{0\}$ &     &     \\ \hline
\end{tabular}
\end{table*}

\subsubsection{Case (2): Input PDF is only partially known through marginals and correlation matrix}
\label{sec:ch4-case2}
Accurate UQ relies on a well-defined joint distribution of input random variables. However, in many practical engineering applications, obtaining such a complete joint distribution is challenging due to the limited availability of statistical data (\citealt{beer2013reliability,kumar2023bayesian}). Often, the available information is only restricted to the marginal distributions and the covariance matrix (or pairwise correlations) (\citealt{der1986structural,li2012uncertainty,dutfoy2009practical}). Deriving the full joint distribution solely from these measures is inherently indeterminate, as the correlation matrix alone does not uniquely specify the complex dependency structure among variables (\citealt{phoon2004simulation,kazianka2011bayesian}).

Under such limited-information settings, the Nataf transformation has been widely adopted in conjunction with the PEM framework (\citealt{chen2015correlated, ahmadabadi2015assessment}). As discussed earlier, the Nataf transformation corresponds to a Gaussian-copula representation of the dependence structure. The prescribed pairwise linear correlation coefficients $\rho_{ij}$ in the physical space are defined as:
\begin{equation}
    \label{eq:def_correlation}
    \rho_{ij}=\int_{-\infty}^{\infty} \int_{-\infty}^{\infty}\left(\frac{x_i-\mu_i}{\sigma_i}\right)\left(\frac{x_j-\mu_j}{\sigma_j}\right) f\left(x_i, x_j\right) d x_i d x_j
\end{equation}
where $\mu_i$ and $\sigma_i$ are the mean and standard deviation of $x_i$, and $f(x_i, x_j)$ represents the bivariate marginal PDF. Nevertheless, the available marginal and correlation information alone does not imply that the underlying dependence structure is Gaussian, and alternative dependence models may be consistent with the same prescribed information. Therefore, the choice of dependence model constitutes an additional modeling assumption under incomplete probability information. Accordingly, this study investigates the sensitivity of QPEM moment estimates to this modeling assumption by considering flexible vine-copula constructions as alternatives to the Gaussian copula.

To construct alternative joint distributions consistent with the available information, a vine structure and corresponding pair-copula families must first be specified. Specifically, we consider two principal classes of regular vines, namely C-vines and D-vines. This choice is motivated by the lack of definitive evidence favoring one over the other when the true structure is unknown. Once a vine structure is selected, a specific pair-copula family is assigned to each edge. For an $n$-dimensional input vector, there are $n!/2$ possible C-vines and $n!/2$ D-vines, each representing a different configuration of pairwise dependencies.

The main challenge lies in inferring the parameters of these pair-copulas using the available correlation matrix. The joint density $f(x_i, x_j)$ in Eq.~\eqref{eq:def_correlation} can be expressed through conditional probabilities as:
\begin{equation}
    \label{eq:bivariate_pdf_conditioning}
    f\left(x_i, x_j\right)=\frac{f(\mathbf{v})}{f(\mathbf{v}_{\overline{i}\overline{j}} \mid x_i, x_j)}
\end{equation}
where both the numerator and denominator can be expressed by products of pair-copulas in a vine structure. Unfortunately, Eq.~\eqref{eq:def_correlation} does not allow for a direct analytical calculation of higher-order pair-copula parameters from $\rho_{ij}$, except for the first level. This is because pair-copulas at higher levels are interdependent, and a direct analytical mapping from the prescribed linear correlations to the corresponding pair-copula parameters is generally unavailable (\citealt{wang2017towards}).

Vine trees, however, provide a structured way to establish a one-to-one mapping between pair-copula parameters and linear correlation coefficients (\citealt{bedford2002vine}). In this framework, each pair-copula parameter corresponds directly to a specific entry in the correlation matrix. Leveraging the hierarchical structure of vine copulas discussed earlier, we adopt a sequential search strategy to infer the pair-copula parameters based solely on the correlation matrix (\citealt{wang2018system}). For a fixed vine and pair-copula specification, the induced correlation is typically monotonic with respect to the corresponding copula parameter over its admissible range. Therefore, when a feasible solution exists, the target parameter can be efficiently identified using a standard one-dimensional root-finding procedure. In this study, a bisection-based calibration procedure is employed. Importantly, this calibration procedure involves only sampling and transformation of the input probabilistic model and does not require evaluations of the computational response model. The procedure involves the following steps:
\begin{itemize}
    \item[(a)] \textbf{Initialize Model:} Select a vine structure (C-vine or D-vine) and assign candidate pair-copula families to each edge. Initialize the pair-copula parameters within their feasible ranges. 
    \item[(b)] \textbf{Transform Sigma Points:} Generate QPEM sigma points in the independent standard normal space, then apply the inverse Rosenblatt transformation based on the current vine copula model. This step involves only the probabilistic transformation of the sigma points and does not require evaluations of the computational model.
    \item[(c)] \textbf{Compute Correlation:} Calculate the linear correlation coefficients $\hat{\rho}_{ij}$ $(i>j)$ from the transformed sigma points.
    \item[(d)] \textbf{Update Parameters:} If the target correlation is attainable under the current pair-copula specification, iteratively update the copula parameters using the bisection method until the estimated correlation $\hat{\rho}_{ij}$ matches the target $\rho_{ij}$ within a specified tolerance.
\end{itemize}
Upon convergence of the sequential search, the resulting pair-copula parameters define a vine model that is consistent with the prescribed correlation matrix within the specified numerical tolerance. The vine-copula framework thereby allows alternative joint dependence structures to be constructed while maintaining consistency with the prescribed marginal distributions and correlation matrix. This provides a means to investigate the sensitivity of the resulting QPEM moment estimates to the dependence-model assumption under incomplete probability information.

\bmsection{Numerical examples}
\label{sec:numerical_examples}

In this section, representative numerical examples are presented to analyze the performance and accuracy of the proposed copula-based QPEM framework. To validate the efficacy of the proposed method, the results are benchmarked against Monte Carlo Simulation (MC) with $10^6$ samples, which serves as the reference solution. In order to compare with the copula-based QPEM with $2n^2+1$ needed samples, variance reduction sampling techniques are also employed. In particular, Latin Hypercube Sampling (LHS) and Quasi-Monte Carlo simulation (QMC) with $2n^2+1$ samples, for the same computational cost are used, together with the Smolyak Gauss-Hermite sparse grid quadrature (SGH3) with $2n^2+2n+1$ required samples, and one of the most well-known and popular PEMs, Hong's PEM (HPEM), with $2n+1$ needed samples. For all numerical examples, the QPEM parameters are set to $r=3$, $\zeta=-8$, and $\xi=60$, consistent with the formulation described earlier.

\subsection{Case (1): Input joint distribution inferred from raw data}
\label{sec:ex1_intro}
This section investigates the application of our proposed framework under Case (1). As defined previously, this scenario addresses the situation where the input joint probability distribution is empirically inferred directly from the available multivariate dataset. In this context, the vine copula structure and its pair-copula families are not assumed a priori, but rather are learned from the data using the sequential heuristic algorithm described in Table~\ref{table:seq_method}.

To evaluate the proposed framework under a controlled finite-sample setting, a synthetic data strategy is employed. For the subsequent numerical illustrations, a reference joint distribution (ground truth) is first established to generate a dataset comprising $N=500$ realizations. This sample size is adopted as a representative finite-data setting for the problem dimensions considered in this study and is not intended to imply a universal sample-size requirement for vine-copula inference. This dataset is then treated as the available information for inferring the vine copula model, which is subsequently used within the QPEM framework to estimate the probabilistic output moments.

To distinguish the effect of finite-sample input-distribution inference from the numerical error associated with moment propagation, two MC reference solutions are considered. The first, denoted as \textit{MC-GT}, is obtained using $10^6$ samples generated directly from the prescribed ground-truth input distribution and provides an end-to-end reference for the complete inference-and-propagation procedure. The second, denoted as \textit{MC-Fitted}, is obtained using $10^6$ samples generated from the probabilistic input model inferred from the $N=500$ dataset. MC-Fitted therefore serves as the reference for assessing the numerical accuracy of QPEM and the competing moment-estimation methods under the same inferred input model. The discrepancy between MC-GT and MC-Fitted reflects the effect of finite-sample vine-copula inference. For compact presentation, the relative errors for all four Case (1) numerical examples are consolidated in Figure~\ref{fig:case1_error}, with panels (a)-(d) corresponding to Examples 1-4, respectively.

\subsubsection{Example 1: Rock slope safety}
\label{sec:ex1-rockslope}
\begin{figure}[!t]
    \centering
    \includegraphics[width=0.6\textwidth]{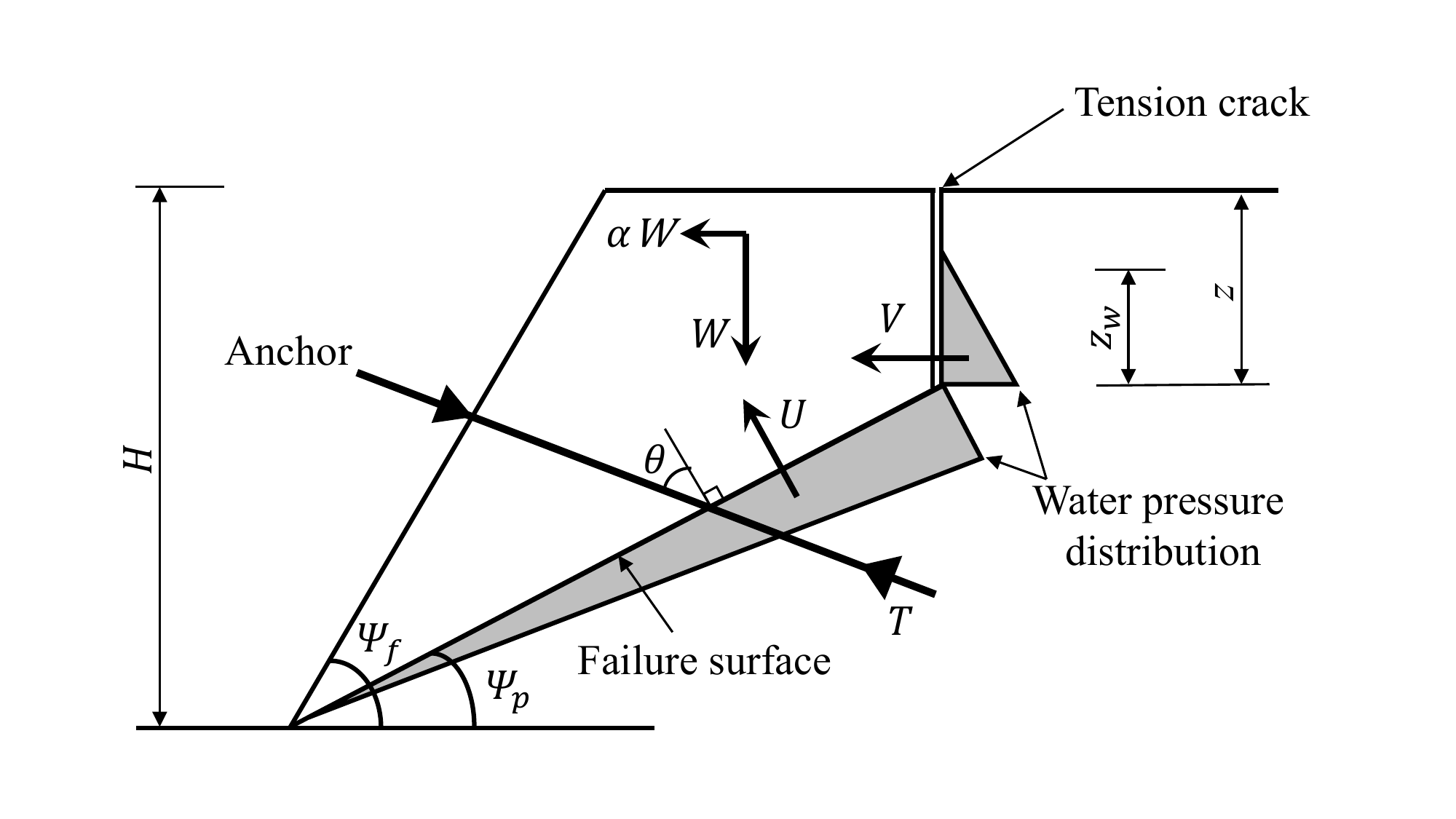}
    \caption{Schematic representation of a rock slope (\citealt{hoek2007practical}).}
\label{fig:ch4-ex1-rockslope}
\end{figure}

This example considers a probabilistic stability analysis of a rock slope with a water-filled tension crack, as illustrated in Figure~\ref{fig:ch4-ex1-rockslope}. The factor of safety ($F_s$) in this case is given by \cite{hoek2007practical}:
\begin{equation}
\label{eq:model_rock_slope}
    F_s(\boldsymbol{X})=\frac{h A+N^{\prime} \tan \phi}{W\left(\sin \psi_p+\alpha \cos \psi_p\right)+V \cos \psi_p-T \sin \psi_a}
\end{equation}
where\\
$\begin{array}{ll}
    \qquad A  &=(H-z) / \sin \psi_p \\
    \qquad z  &=H\left(1-\sqrt{\cot \psi_f \tan \psi_p}\right) \\
    \qquad N^{\prime} &=W\left(\cos \psi_p-\alpha \sin \psi_p\right)-U-V \sin \psi_p+T \cos \theta \\
    \qquad W           &=0.5 \gamma H^2\{[1-\left({z}/{H}\right)^2] \cot \psi_p-\cot \psi_f\} \\
    \qquad U           &=0.5 \gamma_w z_w A \\
    \qquad V           &=0.5 \gamma_w z_w^2 \\
    \qquad r           &=z_w / z 
\end{array}$\\
in which $\gamma$=unit weight of rock=$2.6\times10^4$$N/m^3$, $\gamma_w$=unit weight of water=$1.0\times10^4$$N/m^3$, $\psi_f$=angle of slope=$50^{\circ}$, $\psi_p$=angle of failure surface=$35^{\circ}$, $T$=force applied by anchor system=$0$$N$, $\theta$=inclination of anchor=$0^{\circ}$ and $H$=height of the overall slope=60$m$.

The problem involves a five-dimensional stochastic space ($n$=5) comprising the cohesion $h$ and the friction angle $\phi$ of the failure surface, the tension crack depth $z$, water-depth ratio $r=z_w/z$, and the horizontal earthquake acceleration $\alpha$. The water height $z_w$ is consequently determined as $z_w=rz$. The statistical characteristics of their marginal distributions are summarized in Table~\ref{table:rock_slope_pdf}.

\begin{table}[!b]
\caption{Marginals of random variables for the rock slope example}
\label{table:rock_slope_pdf}
\begin{tabular*}{\columnwidth}[t]{p{0.2\columnwidth} p{0.4\columnwidth} p{0.15\columnwidth} p{0.15\columnwidth}}
\toprule
\textbf{Variable} & \textbf{Distribution} & \textbf{Mean}                   & \textbf{CoV}\\ 
\midrule
$h$ $[kPa]$       & Lognormal    & 140                &   0.10  \\
$\phi$ $[^{\circ}]$       & Lognormal    & 30                   & 0.10            \\
$z$ $[m]$       & Lognormal    & 14  & 0.15                 \\
$r$        & Truncated Exponential [0, 1] & 0.3                & 0.10     \\
$\alpha$ $[g]$    & Truncated Exponential [0, 0.15] & 0.02                & 0.10  \\
\bottomrule
\end{tabular*}
\end{table}

For the establishment of a reference solution, a specific dependence structure based on a D-vine is assumed here, as follows:
\begin{equation}
\label{eq:copula_ex1}
    \begin{split}
        c_{h,\phi, z, r, \alpha}&= c_{h \phi} \cdot c_{z r} \cdot \underbrace{c_{\phi z} \cdot c_{r \alpha}\cdot c_{hz \mid \phi} \cdot c_{\phi r \mid z} \cdot c_{z \alpha \mid r} \cdot c_{hr \mid \phi z} \cdot c_{\phi \alpha \mid z r} \cdot c_{h \alpha \mid \phi z r}}_{= 1} 
    \end{split}
\end{equation}
where the pairs $h-\phi$ and $r-z$ are correlated via a Frank copula rotated by 90$^{\circ}$ ($\theta=2$) and a t-copula ($\rho=-0.6$, $\nu=10$), respectively. No additional variable dependencies ($c$=1) are assumed in this example.


Based on this reference model and the sample size of 500 realizations, the vine copula specification is inferred using the sequential method discussed in Table~\ref{table:seq_method}, and the first four moments of the response function in Eq.~\eqref{eq:model_rock_slope} are estimated using the proposed QPEM. Table~\ref{table:ch4-ex1-rockslope-result} reports two MC reference solutions: MC-GT, obtained directly from the prescribed ground-truth model in Eq.~\eqref{eq:copula_ex1}, and MC-Fitted, obtained from the probabilistic model inferred from the 500 realizations. The two MC solutions remain in close agreement, with discrepancies of approximately 4\% or less across all four response moments, indicating that the finite-sample inferred model provides a reasonable representation of the prescribed ground-truth distribution for this example.


\begin{table}[!tb]
	\caption{Moment estimations for the rock slope example ($n=5$)}
	\label{table:ch4-ex1-rockslope-result} 
	\begin{tabular*}{\textwidth}[t]{p{0.13 \textwidth} p{0.17\textwidth} p{0.13\textwidth} p{0.13\textwidth} p{0.13\textwidth} p{0.13\textwidth}}
		\toprule
		\textbf{Method}   & \centering \textbf{No. of points}& \textbf{Mean}& \textbf{STD} & \textbf{Skewness} & \textbf{Kurtosis}\\\midrule
		MC-GT  & \centering $10^6$ & 1.5367 & 0.1223		& 0.1381 & 3.3326\\ 
        MC-Fitted  & \centering $10^6$ & 1.5407 & 0.1274		& 0.1436 & 3.4628\\ 
	LHS & \centering51 & 1.5429       & 0.1357        & 0.8177        & 4.0762 \\
	QMC  & \centering51 & 1.5389       & 0.1025        & 0.3359        & 2.5612 \\
	SGH3  & \centering61   & 1.5404      & 0.1250        & 0.1123        & 2.6267\\
	HPEM  & \centering11 & 1.5417       & 0.1347        & 0.4462        & 2.1086 \\
	QPEM   & \centering51  & 1.5406       & 0.1257        & 0.1360        & 3.3235 \\ \bottomrule
	\end{tabular*}
\end{table}


\begin{figure}[!t]
    \centering
    \includegraphics[width=\textwidth]{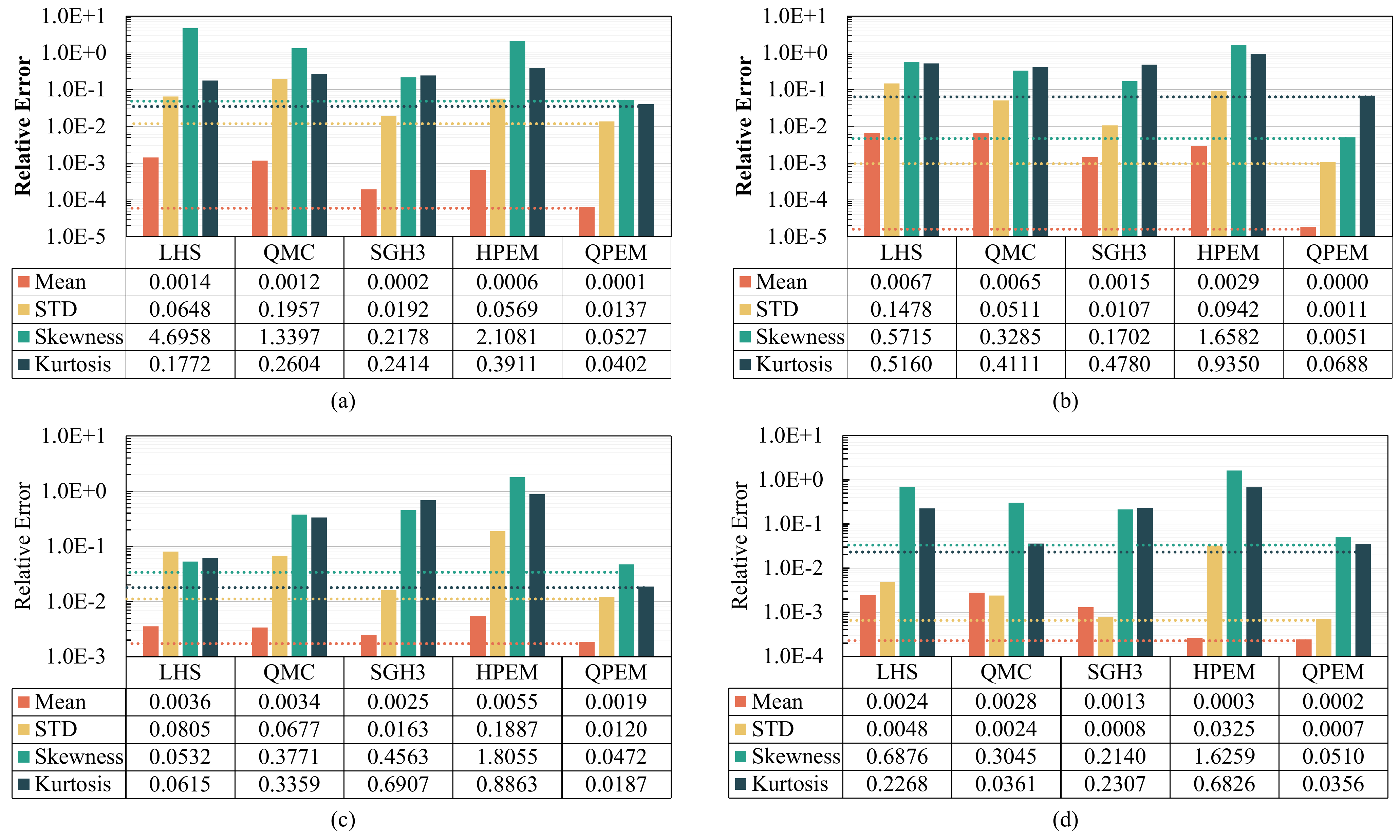}
    \caption{Relative errors for the Case (1) numerical examples: (a) rock slope; (b) rock tunnel excavation; (c) horizontal truss; and (d) space truss.}
\label{fig:case1_error}
\end{figure}

\indent The relative errors in Figure~\ref{fig:case1_error}(a) are evaluated with respect to MC-Fitted to separate the numerical moment-propagation error from the discrepancy introduced during input-distribution inference. Generally, the first two output moments are estimated reasonably well by most methods. Regarding the higher-order moments (skewness and kurtosis), however, all competing methods exhibit considerably larger relative errors than the QPEM. In particular, LHS and QMC, based on the same 51 model evaluations as the QPEM ($2\times5^2+1$), provide substantially less accurate higher-order moment estimates, while SGH3 and HPEM also do not achieve comparable accuracy. In contrast, the QPEM yields relative errors of approximately 0.01\% and 1.3\% for the mean and standard deviation, respectively, and approximately 5.3\% and 4.0\% for skewness and kurtosis. Overall, these results demonstrate that the proposed QPEM provides accurate estimation of the first four response moments while maintaining a favorable balance between computational cost and accuracy.

\subsubsection{Example 2: Rock tunnel excavation}
\label{sec:ch4-rock-tunnel-excavation}
\begin{figure}[!t]
    \centering
    \includegraphics[width=\textwidth]{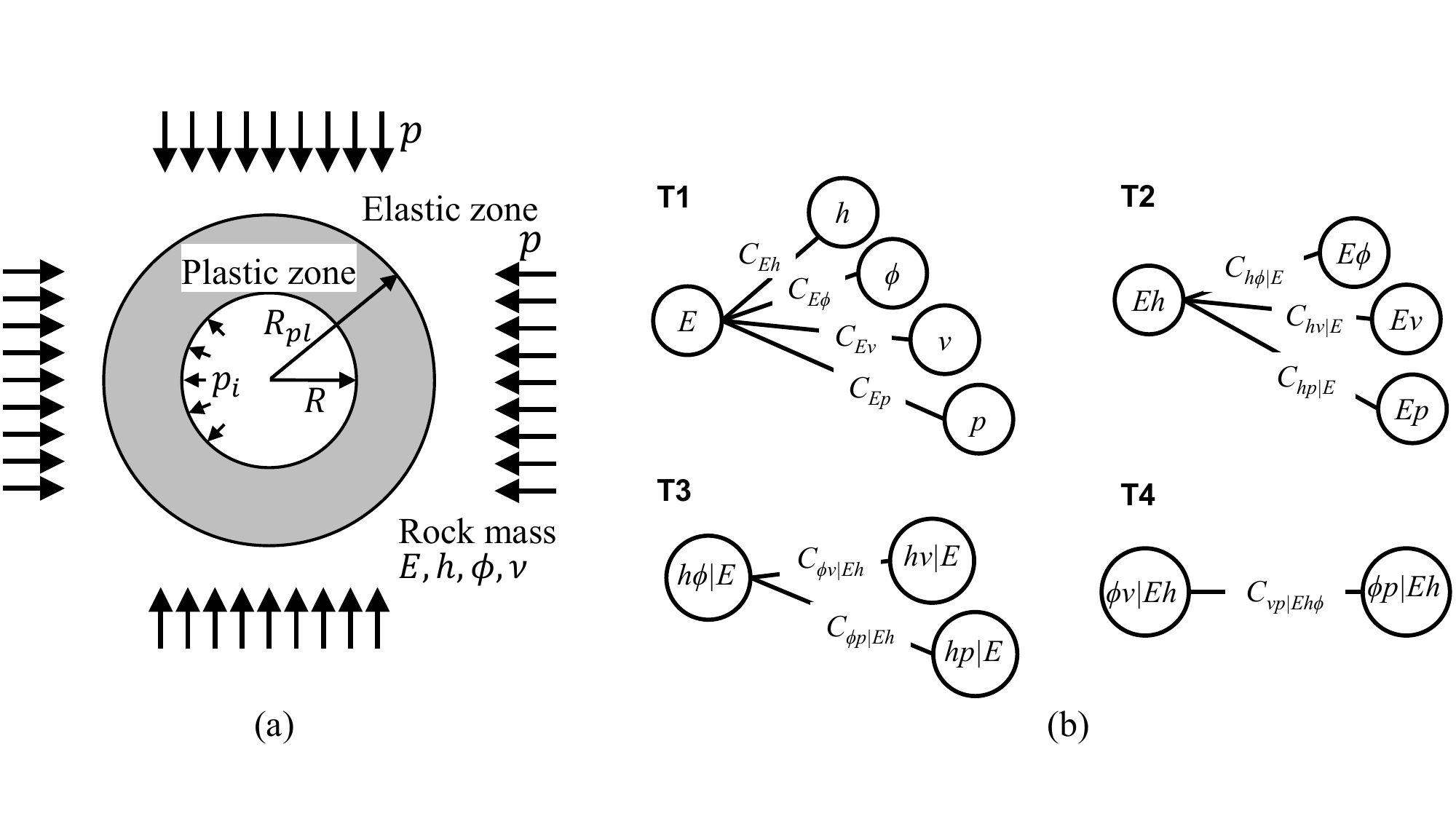}
    \caption{Rock tunnel excavation example: (a) schematic representation of the circular tunnel excavation in a rock mass; and (b) C-vine dependence structure adopted for the stochastic input variables.}
\label{fig:ch4-ex2-rocktunnel-drawing}
\end{figure}
This numerical example addresses a rock tunnel excavation, adapted from \cite{li2010reliability}. As schematically shown in Figure~\ref{fig:ch4-ex2-rocktunnel-drawing}(a), the excavation process induces significant stress redistribution and displacements within the surrounding rock mass. The analytical solution for the induced inward displacement, derived based on the Mohr-Coulomb failure criterion, is governed by the relationship between the support pressure $p_i$ and the critical pressure $p_{cr}$. The mechanical behavior of the tunnel exhibits two distinct regimes depending on the magnitude of the uniform support pressure $p_i$ relative to the critical pressure $p_{cr}$. The critical pressure is defined as $p_{cr}=({2 p-\sigma_c})/({k+1})$, where the auxiliary parameters are given by $k=({1+\sin \phi})/({1-\sin \phi})$ and $\sigma_c={h(k-1)}/{\tan \phi}$. Here, $h$ denotes the cohesion and $\phi$ represents the friction angle.

When $p_i < p_{cr}$, a plastic zone develops around the excavation. The radius of this plastic zone, $R_{pl}$, is expressed as:
\begin{equation}
    {R_{p l}}=R\left[\frac{2(p+s)}{(k+1)\left(p_i+s\right)}\right]^{1 /(k-1)}
\end{equation}
where $R$ is the radius of the circular tunnel, $p$ is the hydrostatic pressure, and $s={\sigma_c}/({k-1})$. In this plastic regime, the inward displacement of the tunnel lining, $u_{ip}$, is calculated as:
\begin{equation}
    u_{i p}=R\left[\frac{(1+v)}{E}\right]\left[2(1-v)\left(p-p_{c r}\right)\left(\frac{R_{p l}}{R}\right)^2-(1-2 v)\left(p-p_i\right)\right]
\end{equation}
where $E$ is the Young's modulus and $v$ is the Poisson's ratio. Conversely, if the support pressure is sufficient to prevent plastic yielding (i.e., $p_i \geqslant p_{c r}$), the displacement follows the elastic solution:
\begin{equation}
    u_{i p}=R\left[\frac{(1+v)}{E}\right]\left(p-p_i\right)
\end{equation}

In this study, the uniform support pressure $p_i$ is treated as a deterministic parameter with a value of 0.868 MPa. It is important to note that since $p_{cr}$ is instead a function of random variables, the system exhibits a stochastic regime-switching behavior. This implies that the response function $g(\mathbf{X})$ involves a gradient discontinuity at the transition point between the elastic and plastic states. Such non-smoothness typically poses significant challenges for conventional numerical integration methods. The response function is defined based on a serviceability criterion, where the permissible inward displacement is restricted to $1\%$ of the tunnel radius. Accordingly, the response model is formulated as:
\begin{equation}
    \label{eq:ch4-tunnel-excavation-function}
    g(\boldsymbol{X})=0.01-\frac{u_{i p}}{R}
\end{equation}

\begin{table}[t]
\caption{Marginals of random variables for the tunnel excavation example.}
\label{table:ch4-ex2-rocktunnel-variable}
\begin{tabular*}{\columnwidth}[t]{p{0.25\columnwidth} p{0.25\columnwidth} p{0.25\columnwidth} p{0.25\columnwidth}}
\toprule
\textbf{Variable} & \textbf{Distribution} & \textbf{Mean}   & \textbf{CoV} \\ 
\midrule
$E$ $[MPa]$       & Lognormal    & 373 & 0.13  \\
$h$ $[MPa]$       & Lognormal    & 0.23  & 0.10  \\
$\phi$ $[^{\circ}]$       & Lognormal    & 22.85  & 0.06    \\
$v$        & Beta [0.2, 0.4] & 0.3                & 0.17      \\
$p$ $[MPa]$    & Beta [1.5, 3]& 2                & 0.15      \\ 
\bottomrule
\multicolumn{4}{l}{\begin{tabular}{@{}l@{}}  \end{tabular}} 
\end{tabular*}
\end{table}


The probabilistic model involves five random variables ($n=5$), the Young's modulus $E$, cohesion $h$, friction angle $\phi$, Poisson's ratio $v$, and hydrostatic pressure $p$. The marginal distributions for these variables are summarized in Table~\ref{table:ch4-ex2-rocktunnel-variable}. To model the dependence structure, a C-vine topology is utilized as illustrated in Figure~\ref{fig:ch4-ex2-rocktunnel-drawing}(b). Specifically, Young's modulus $E$ is designated as the central root node in the first tree $T_1$, a choice justified by its governing role in the overall displacement magnitude. A diverse set of pair-copulas is assigned to the edges to introduce complex dependency features. The first tree includes $C_{E h}$ and $C_{E \phi}$ modeled by Clayton copulas ($\theta=1.2$), $C_{E p}$ by a Clayton copula ($\theta=0.6$), and $C_{E v}$ by a Frank copula rotated by 90$^{\circ}$ ($\theta=1.5$). Higher-order dependencies include $C_{h \phi \mid E}$ modeled as a Frank copula rotated by 90$^{\circ}$ ($\theta=2.0$), while the remaining pairs ($C_{h v \mid E}$, $C_{h p \mid E}$, $C_{\phi v \mid E h}$, $C_{\phi p \mid E h}$, and $C_{v p \mid E h \phi}$) are assumed to be Gaussian with $\rho=0$. Consistent with the methodology applied in the previous example, the vine copula specification is inferred using the sequential method outlined in Table~\ref{table:seq_method} using the 500 available random samples generated from the defined vine copula model.


The first four moments of the response function are estimated using all considered methods. Table~\ref{table:ch4-ex1-rockstunnel-result} reports both MC-GT, obtained from the prescribed ground-truth input model, and MC-Fitted, obtained from the probabilistic model inferred from the 500 available realizations. The relative errors in Figure~\ref{fig:case1_error}(b) are evaluated with respect to MC-Fitted so that the numerical moment-propagation accuracy can be distinguished from the discrepancy introduced during finite-sample input-distribution inference.

The effect of finite-sample input-model inference is more noticeable for the higher-order response moments in this example. In particular, the discrepancies between MC-GT and MC-Fitted are approximately 6.6\% for skewness and 10.2\% for kurtosis, compared with approximately 4.0\% and 3.9\%, respectively, in the previous example. Nevertheless, QPEM closely reproduces the MC-Fitted mean, standard deviation, and skewness, with relative errors below approximately 0.6\%. The kurtosis is more challenging, with a QPEM relative error of approximately 6.9\%, but this error remains substantially smaller than those of the other methods. Interestingly, the kurtosis of QPEM is closer to MC-GT than to MC-Fitted because the input-model inference discrepancy and the subsequent QPEM approximation error act in opposite directions. This observation further motivates separating input-model inference error from numerical moment-propagation error when assessing the overall framework. Overall, QPEM provides the most accurate representation of the first four moments of the inferred probabilistic model among the considered methods, including for this non-smooth response characterized by stochastic elastic-plastic regime switching.

\begin{table}[!t]
    \caption{Moment estimations for the tunnel excavation example ($n=5$)}
    \label{table:ch4-ex1-rockstunnel-result} 
    \begin{tabular*}{\textwidth}[t]{p{0.13 \textwidth} p{0.17\textwidth} p{0.13\textwidth} p{0.13\textwidth} p{0.13\textwidth} p{0.13\textwidth}}
		\toprule
		\textbf{Method}   & \centering \textbf{No. of points}& \textbf{Mean}& \textbf{STD} & \textbf{Skewness} & \textbf{Kurtosis}\\\midrule
		MC-GT & \centering  $10^6$& 0.0059 & 0.0009		& -1.0389 & 4.8973\\ 
        MC-Fitted  & \centering $10^6$ & 0.0058 & 0.0010		& -1.1070 & 5.3954\\ 

		LHS  & \centering  51 & 0.0059        & 0.0008        & -0.4744        & 2.6117 \\
		QMC  & \centering  51& 0.0059        & 0.0009        & -0.7434        & 3.1781 \\
		SGH3    & \centering  61 & 0.0059       & 0.0010        & -0.9186        & 2.8172\\
		HPEM  & \centering  11 & 0.0058        & 0.0009       & 0.7286        & 0.3506 \\
		QPEM    & \centering  51  & 0.0058        & 0.0010        & -1.1126        & 5.0251 \\
        \bottomrule
	\end{tabular*}
\end{table}

\subsubsection{Example 3: Maximum vertical displacement of horizontal truss}
\label{sec:ch4-horizontal-truss}
\begin{figure}[!t]
    \centering
    \includegraphics[trim=0in 1.0in 0in 1.5in, clip=true, width=\textwidth]{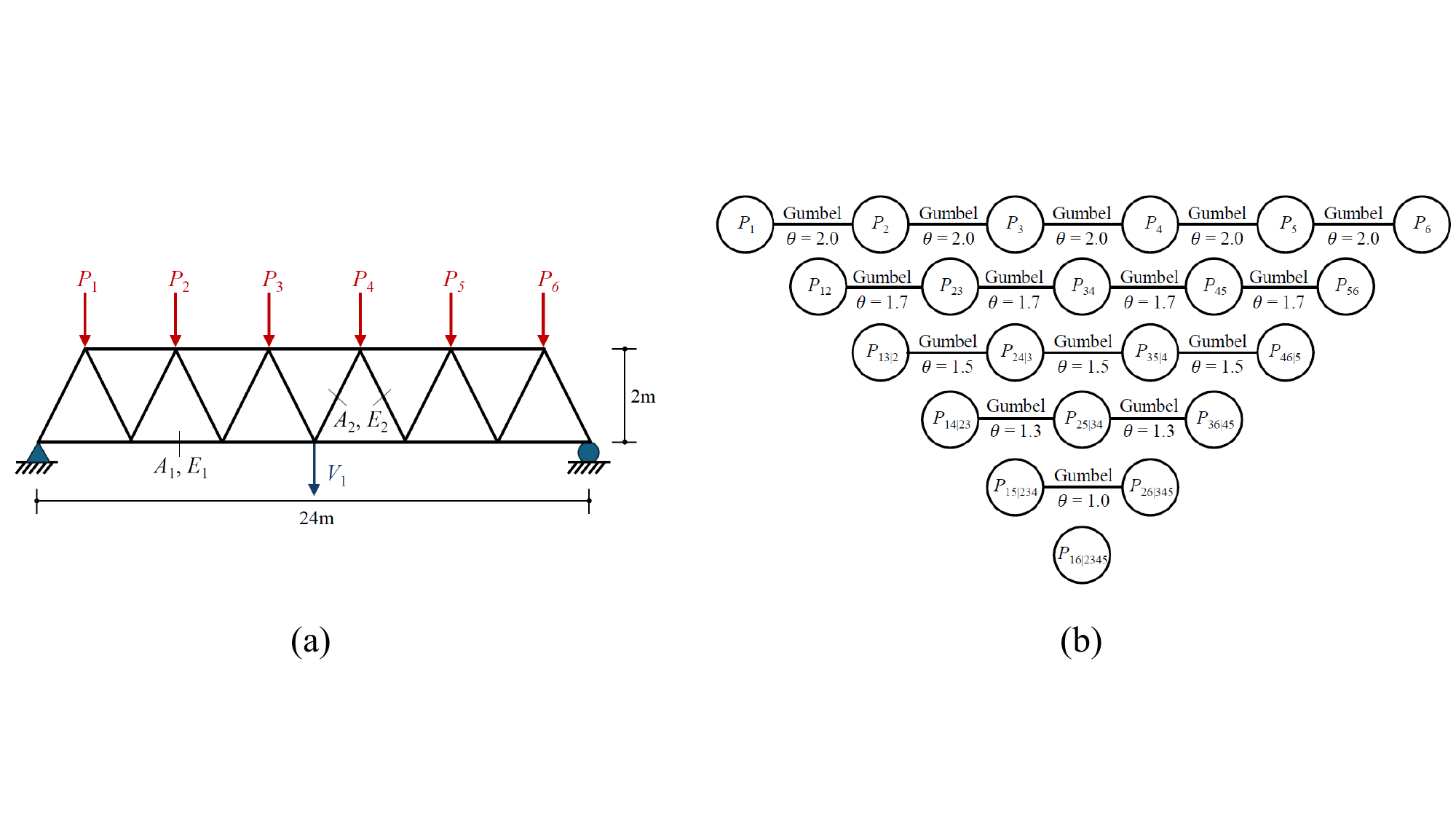}
    \caption{Horizontal truss example: (a) 23-member plane truss (\citealt{blatman2010}); and (b) D-vine dependence structure for the loads $P_1$ through $P_6$}
\label{fig:ch4-horizontal-truss-drawing}
\end{figure}

This numerical example investigates the probabilistic response of the simply supported truss structure adapted from \cite{blatman2010} (Figure~\ref{fig:ch4-horizontal-truss-drawing}(a)). The structure comprises 23 bar elements, divided into 11 horizontal bars and 12 inclined bars, connected by 13 nodes. The system is subjected to vertical point loads, $P_i$, applied to the nodes of the upper horizontal chord. The quantity of interest is the maximum vertical displacement of the structure, denoted as $V_1$. This response is computed using a two-dimensional linear elastic finite element model, defined as:

\begin{equation}
    y=\mathcal{M(\boldsymbol{X})}=V_1(\boldsymbol{X})
\end{equation}

The probabilistic model of the system involves 10 random variables ($n=10$), as summarized in Table~\ref{table:ch4-horizontal-truss-variable}. These variables include the Young’s modulus ($E_1, E_2$) and the cross-sectional areas ($A_1, A_2$) for the horizontal and inclined bar groups, respectively, as well as the six vertical loads ($P_1, \ldots, P_6$).

\begin{table}[!bt]
\caption{\label{table:ch4-horizontal-truss-variable}Random variables in the horizontal truss example}
\centering
\begin{tabular*}{\columnwidth}[t]{p{0.16\columnwidth} p{0.32\columnwidth} p{0.16\columnwidth} p{0.16\columnwidth} p{0.20\columnwidth}}
\toprule
\textbf{Variable} &\textbf{Description}& \textbf{Distribution} & \textbf{Mean} &\textbf{CoV}\\ 
\midrule
$E_1$, $E_2$ [$Pa$]&Young's modulus& Lognormal  & $2.10 \times 10^{11}$ & 0.17 \\
$A_1$ [$m^2$]&Cross section of horizontal bars& Lognormal  & $2 \times 10^{-3}$ & 0.17\\
$A_2$ [$m^2$]&Cross section of inclined bars& Lognormal  & $1 \times 10^{-3}$ & 0.17 \\
$P_1,\ldots,P_6$ [$N$]&Vertical Loads& Gumbel  & $5 \times 10^{4}$ & 0.20 \\ \bottomrule
\end{tabular*}
\end{table}

To characterize the dependence structure within the ten-dimensional stochastic input, a block-independence assumption is adopted (\citealt{blatman2010}). Specifically, it is assumed that the group of material property variables ($E_1, E_2, A_1, A_2$) is statistically independent from the group of load variables ($P_1, \ldots, P_6$). Regarding the internal structure of each block, the material properties are considered mutually independent ($c_{E_1, E_2, A_1, A_2}=1$), whereas the load variables exhibit significant dependence, which is modeled using the D-vine structure illustrated in Figure~\ref{fig:ch4-horizontal-truss-drawing}(b).
The joint copula density is thus expressed as follows:
\begin{equation}
    \begin{split}
        \label{eq:ch4-horizontal-truss-exact-vine-copulas}
        c_{E_1, E_2, A_1, A_2, P_1, \ldots, P_6}&=\underbrace{c_{E_1, E_2, A_1, A_2}}_{=1} \cdot c_{P_1, \ldots, P_6}: \text{block independence}\\
        &=c_{P_1, \ldots, P_6}\\
    &= \underbrace{c_{12} \cdot c_{23} \cdot c_{34} \cdot c_{45} \cdot c_{56}}_{=\text{Gumbel(2.0)}} \cdot \underbrace{c_{13 \mid 2} \cdot c_{24 \mid 3} \cdot c_{35 \mid 4} \cdot c_{46 \mid 5}}_{=\text{Gumbel(1.7)}} \cdot \\ 
    & \qquad \underbrace{c_{14 \mid 23} \cdot c_{25 \mid 34} \cdot c_{36 \mid 45}}_{=\text{Gumbel(1.5)}}\cdot \underbrace{c_{15 \mid 234} \cdot c_{26 \mid 345}}_{=\text{Gumbel(1.3)}} \cdot \underbrace{c_{16 \mid 2345}}_{=\text{Gumbel(1.0)}}
    \end{split}
\end{equation}

\begin{table}[!t]
	\caption{\label{table:ch4-horizontal-truss-results}Moment estimations for the maximum vertical displacement of the horizontal truss structure ($n=10$)}
	\begin{tabular*}{\textwidth}[t]{p{0.15\textwidth} p{0.2\textwidth} p{0.14\textwidth} p{0.12\textwidth} p{0.12\textwidth} p{0.12\textwidth}}
		\toprule
		\textbf{Method}   &\centering \textbf{No. of Points} & \textbf{Mean [mm]}& \textbf{STD [mm]} & \textbf{Skewness}&\textbf{Kurtosis} \\\midrule
		MC-GT    &\centering  $10^6$   & 8.1597 & 2.3359		& 1.1364 & 5.5654 \\ 
        MC-Fitted    &\centering  $10^6$   & 8.1430 & 2.3273		& 1.2110 & 5.4855 \\ 
		LHS     &\centering  201      & 8.1140     & 2.1400       & 1.1466        & 5.8229 \\
		QMC     &\centering  201       & 8.1153      & 2.1697        & 0.7543        & 3.6431 \\
		SGH3    &\centering  221        & 8.1635       & 2.2895        & 0.6584        & 1.6967 \\
		HPEM    &\centering  21        &8.0986      & 1.8883        & -0.9754        & 0.6238 \\
		QPEM ($r=3$)&\centering  201        & 8.1581      & 2.3553        & 1.1539        & 5.3829 \\\bottomrule
	\end{tabular*}
\end{table}


The results are summarized in Table~\ref{table:ch4-horizontal-truss-results}, with the corresponding relative errors shown in Figure~\ref{fig:case1_error}(c). A noticeable discrepancy between MC-GT and MC-Fitted is observed for the skewness, while the other response moments remain relatively close. Despite this finite-sample effect, the QPEM provides accurate estimates of all four moments, with relative errors of approximately 0.2\%, 1.2\%, 4.7\%, and 1.9\% for the mean, standard deviation, skewness, and kurtosis, respectively, with respect to MC-Fitted. In particular, the higher-order moments are estimated considerably more accurately than those obtained from the other methods. Moreover, the QPEM estimates remain within approximately 3.3\% of the MC-GT values for all four moments. These results demonstrate that the QPEM maintains accurate higher-order moment estimation in this ten-dimensional problem using only 201 model evaluations, even when the input dependence model is inferred from a finite dataset.

\subsubsection{Example 4: Maximum displacement of space truss structure}
\begin{figure}[!b]
    \centering
    \includegraphics[width=0.8\textwidth]{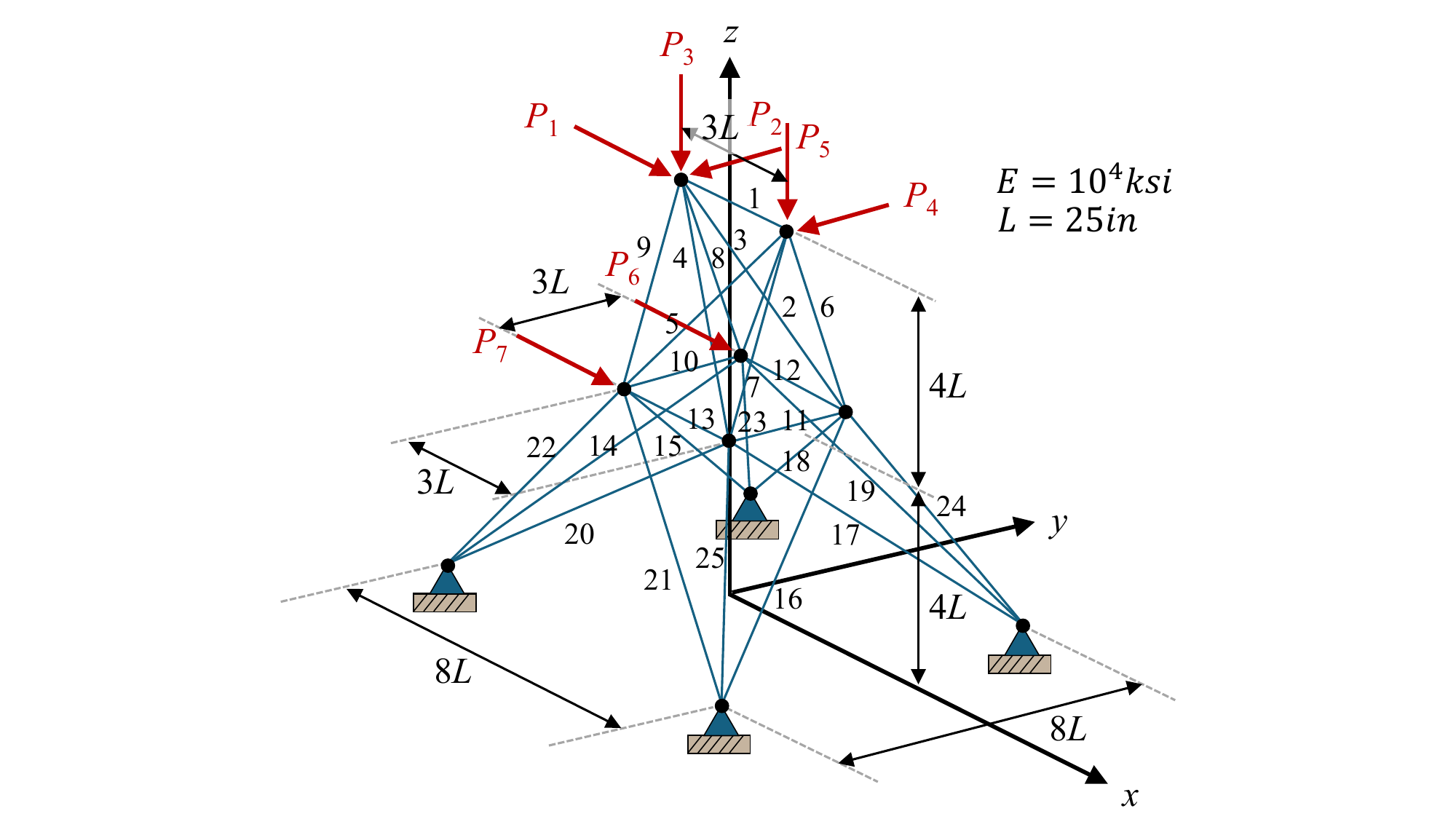}
    \caption{Space truss structure (\citealt{rahami2008sizing}).}
\label{fig:ch4-space-truss-structure}
\end{figure}

This engineering application considers the response of the space truss structure illustrated in Figure~\ref{fig:ch4-space-truss-structure}, as adapted from \citet{rahami2008sizing}. The structural system consists of twenty-five bar elements connected by ten nodes. The truss is subjected to a complex loading scenario comprising five horizontal loads ($P_1, P_2, P_4, P_6, P_7$) and two vertical loads ($P_3, P_5$). The Young's modulus of the material is treated as a deterministic parameter with a value of $10^4$ $ksi$.

The quantity of interest is the maximum nodal displacement at the top of the structure. This response is evaluated using a linear finite element analysis and is formally defined as:
\begin{equation}
    Y=\mathcal{M}(\boldsymbol{X})=\max \left(u_h(\boldsymbol{X}), u_v(\boldsymbol{X})\right)
\end{equation}
where $u_h$ and $u_v$ denote the displacement components in the horizontal and vertical directions, respectively.

The stochastic input vector $\boldsymbol{X}$ comprises 15 random variables ($n=15$), denoted as $\boldsymbol{X}=\left\{X_1, X_2, \ldots, X_{15}\right\}^{T}$. These variables characterize the uncertainties associated with the external loads and the element cross-sectional areas. The probabilistic definitions, including marginal distributions and first two statistical moments, are detailed in Table~\ref{table:ch4-space-truss-structure-variables}. Specifically, variables $X_1$ through $X_7$ represent the applied loads, while variables $X_8$ through $X_{15}$ represent the cross-sectional areas of the structural members, grouped according to symmetry and structural function.

To define the dependence structure among the stochastic inputs, a specific vine copula model is adopted as the exact reference. For the load variables ($X_1-X_7$), the exact dependence is modeled using Gumbel pair-copulas with a parameter $\theta=1.5$. A defining characteristic of this reference model is that the vine structure is truncated after the first tree. This implies that only the unconditional pairwise dependencies are active, while all conditional dependencies in higher-order trees are explicitly defined as independent (i.e., independent copulas). The cross-sectional area variables ($X_8-X_{15}$) are modeled as statistically independent groups.

\begin{table}[!t]
    \caption{Random variables in the space truss structure example}\label{table:ch4-space-truss-structure-variables}
    \begin{tabular*}{\textwidth}[t]{p{0.13\textwidth} p{0.35\textwidth} p{0.15\textwidth} p{0.15\textwidth} p{0.2\textwidth}}
    \toprule
        \centering\textbf{Variable} &\centering \textbf{Description}& \textbf{Distribution} & \centering \textbf{Mean} &  \textbf{CoV}\\ 
    \midrule
        \centering $X_1$ & \centering Horizontal load $P_1$ $[kip]$& \quad Lognormal  & \centering 1 & \;\;0.1 \\
        \centering $X_2$ & \centering Horizontal load $P_2$ $[kip]$& \quad Lognormal  & \centering 12 & \;\;0.1 \\
        \centering $X_3$ & \centering Vertical load $P_3$ $[kip]$& \quad Lognormal  & \centering 10 & \;\;0.1 \\
        \centering $X_4$ & \centering Horizontal load $P_4$ $[kip]$& \quad Lognormal  & \centering 12 & \;\;0.1 \\
        \centering $X_5$ & \centering Vertical load $P_5$ $[kip]$& \quad Lognormal  &\centering  10 & \;\;0.1 \\
        \centering $X_6$ & \centering Horizontal load $P_6$ $[kip]$& \quad Lognormal  & \centering 0.5 & \;\;0.1 \\
        \centering $X_7$ & \centering Horizontal load $P_7$ $[kip]$& \quad Lognormal  &\centering  0.6 & \;\;0.1 \\
        \centering $X_8$ & \centering Cross section $A_1$ $[in^2]$& \quad Lognormal  &\centering  0.4 & \;\;0.1 \\
        \centering $X_9$ & \centering Cross sections $A_2 - A_5$ $[in^2]$& \quad Lognormal  & \centering 0.1 & \;\;0.1 \\
        \centering $X_{10}$ & \centering Cross sections $A_6 - A_9$ $[in^2]$& \quad Lognormal  &\centering  3.4 & \;\;0.1 \\
        \centering $X_{11}$ & \centering Cross sections $A_{10}, A_{11}$ $[in^2]$& \quad Lognormal  & \centering 0.4 & \;\;0.1 \\
        \centering $X_{12}$ & \centering Cross sections $A_{12}, A_{13}$ $[in^2]$& \quad Lognormal  & \centering 1.3 & \;\;0.1 \\
        \centering $X_{13}$ & \centering Cross sections $A_{14}- A_{17}$ $[in^2]$& \quad Lognormal  & \centering 0.9 & \;\;0.1 \\
        \centering $X_{14}$ & \centering Cross sections $A_{18}- A_{21}$ $[in^2]$& \quad Lognormal  & \centering 1.0 & \;\;0.1 \\
        \centering $X_{15}$ & \centering Cross sections $A_{22}- A_{25}$ $[in^2]$& \quad Lognormal  & \centering 3.4 & \;\;0.1 \\ \bottomrule
    \end{tabular*}
\end{table}



Again, the MC-Fitted solution remains close to MC-GT across all four response moments, indicating that the finite-sample inference introduces only a limited discrepancy at the response level. As reported in Table~\ref{table:ch4-space-truss-results} and Figure~\ref{fig:case1_error}(d), some of the competing methods provide accurate estimates for individual statistics but show considerably larger deviations for others. In particular, SGH3 accurately estimates the standard deviation but exhibits substantial errors in skewness and kurtosis, while QMC provides a comparable kurtosis estimate but a much less accurate skewness estimate. In contrast, the QPEM provides consistently accurate estimates across all four moments, with relative errors of approximately 0.02\%, 0.07\%, 5.1\%, and 3.6\% for the mean, standard deviation, skewness, and kurtosis, respectively. This balanced performance is achieved with only 451 model evaluations.

\begin{table}[!t]
	\caption{\label{table:ch4-space-truss-results}Moment estimations for the maximum displacement of the space truss structure ($n=15$)}
	\begin{tabular*}{\textwidth}[t]{p{0.17\textwidth} p{0.2\textwidth} p{0.12\textwidth} p{0.12\textwidth} p{0.12\textwidth} p{0.12\textwidth}}
		\toprule
		\textbf{Method}   &\centering \textbf{No. of Points} & \textbf{Mean [$in$]}& \textbf{STD [$in$]} & \textbf{Skewness}& \textbf{Kurtosis}\\\midrule
		MC-GT    &\centering  $10^6$   & 0.1137 & 0.0160		& 0.2970 & 3.1795 \\ 
        MC-Fitted    &\centering  $10^6$   & 0.1137 & 0.0162		& 0.3047 & 3.2232 \\ 
		LHS    &\centering  451      & 0.1135      & 0.0163       & 0.5142        & 3.9543 \\
		QMC    &\centering  451       & 0.1134     & 0.0163        & 0.2119       & 3.3395 \\
		SGH3   &\centering  481        & 0.1136       & 0.0162        & 0.2395        & 2.4795 \\
		HPEM   &\centering  31        &0.1137      &0.0157        & -0.1907        & 1.0231 \\
		QPEM ($r=3$)&\centering  451        & 0.1137      & 0.0162        & 0.2892        & 3.3380 \\\bottomrule
	\end{tabular*}
\end{table}


\subsection{Case (2): Only marginals and correlation matrix are known}
In this section, we investigate the performance of the copula-based QPEM under the conditions of Case (2), where the probabilistic input is characterized only partially, specifically through its marginal distributions and the correlation matrix. The numerical examples presented here are counterparts to those analyzed for Case (1), retaining the same definitions for random variables and underlying dependence structures.

The fundamental premise for Case (2) is that the available information is strictly limited to the marginals and the correlation matrix. Unlike the previous case where these properties were inferred from a dataset, here the target marginals and correlations are only given. For the purpose of this study, the prescribed marginal distributions are taken directly from the ground-truth models defined in the corresponding Case (1) examples. The target correlation matrices are evaluated numerically from the corresponding ground-truth input models and are provided in Appendix~\ref{app:target_corr}. These calculations involve only sampling from the prescribed input distributions and do not require evaluations of the computational response models. Leveraging this limited information (i.e., the target marginals and the correlation matrix), the first four moments of the response are computed using the vine copula approach described above.

Unlike Case (1), a separate MC-Fitted reference is not considered because the input dependence model is not inferred from finite data in Case (2). Instead, the assumed dependence model represents an explicit modeling choice under incomplete probability information, and the MC-GT solution from the corresponding Case (1) example is used as a common benchmark for the resulting response predictions. The performance of the copula-based QPEM, which requires $2n^2+1$ samples, is again compared against LHS and QMC, both utilizing $2n^2+1$ samples; SGH3, which demands $2n^2+2n+1$ samples in these examples; and HPEM, which operates with $2n+1$ samples.

\subsubsection{Example 5: Maximum vertical displacement of horizontal truss}
The horizontal truss example introduced in Example 3 is now revisited. In this section, it is assumed that our knowledge of the input distribution is strictly limited to the exact marginal distributions and the correlation matrix.

\begin{table}[!t]
	\caption{\label{table:ch4-horizontal-truss-gauss-results}Moment estimations for the maximum vertical displacement of the horizontal truss structure modeled by Gaussian copula ($n=10$)}
	\begin{tabular*}{\textwidth}[t]{p{0.15\textwidth} p{0.2\textwidth} p{0.14\textwidth} p{0.12\textwidth} p{0.12\textwidth} p{0.12\textwidth}}
		\toprule
		\textbf{Method}   &\centering \textbf{No. of Points} & \textbf{Mean [mm]}& \textbf{STD [mm]} & \textbf{Skewness} & \textbf{Kurtosis}\\\midrule
		MC     &\centering  $10^6$   & 8.1597   & 2.3359		& 1.1364 & 5.5654 \\ 
		LHS    &\centering  201      & 8.1187     & 2.4828       & 1.2969        & 3.9216 \\
		QMC    &\centering  201       & 8.1924      & 2.2020        & 0.8303        & 3.9826 \\
		SGH3   &\centering  221        & 8.1231       & 2.4728        & 0.6809        & 1.8596 \\
		HPEM   &\centering  21        &8.1076      & 2.1761        & -0.2147        & 0.8930 \\
		QPEM ($r=3$)&\centering  201  & 8.1322      & 2.3739        & 1.0326        & 4.6770 \\\bottomrule
	\end{tabular*}
\end{table}
\begin{figure}[!t]
    \centering
    \includegraphics[width=0.6\textwidth]{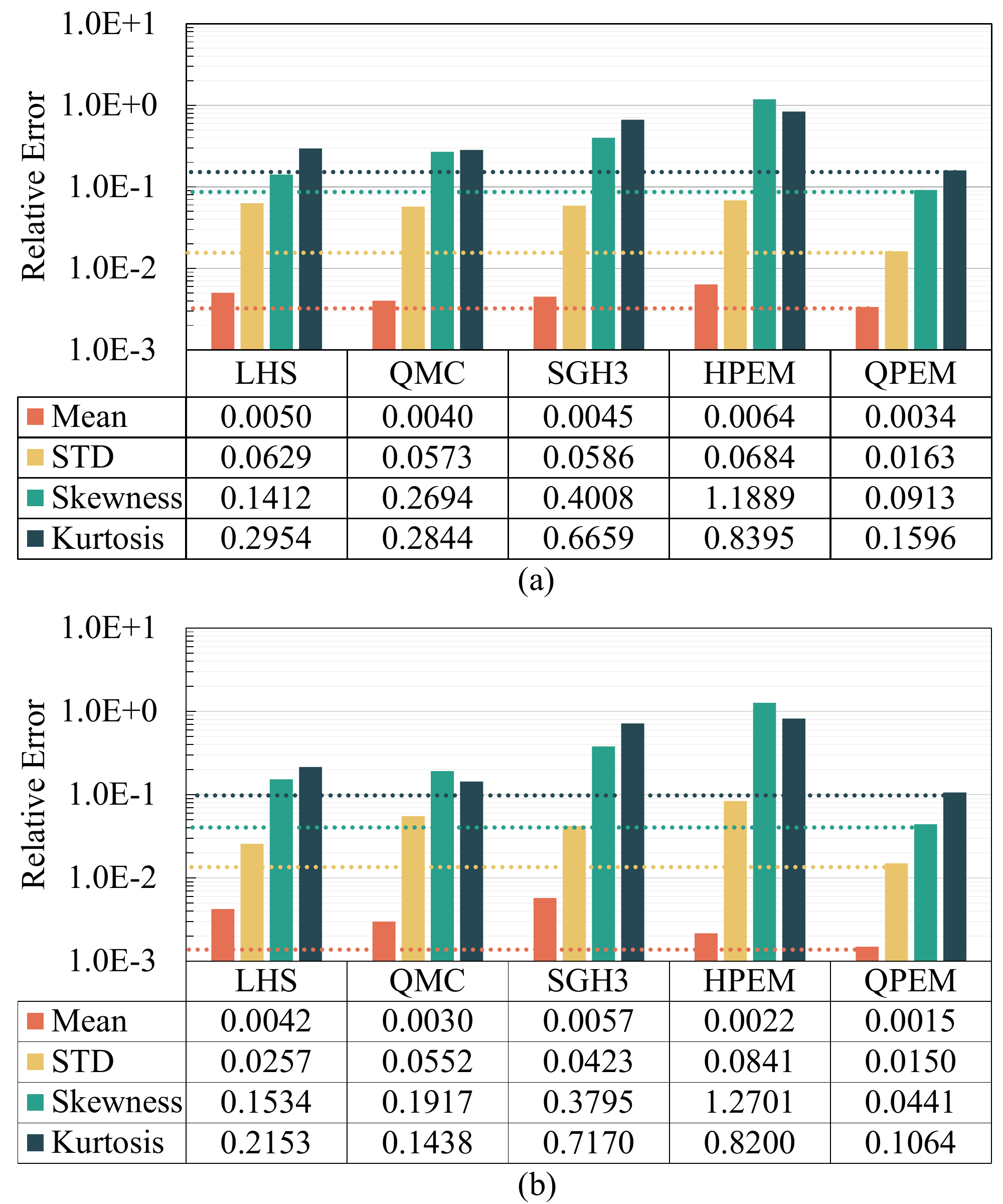}
    \caption{Relative errors for the horizontal truss example under different dependence models: (a) Gaussian copula; and (b) Gumbel copula.}
\label{fig:ch4-horizontal-truss-copula}
\end{figure}



Given the prescribed marginal distributions and correlation matrix, a Gaussian copula is first adopted as a baseline dependence model, corresponding to the conventional Nataf transformation. The Gaussian pair-copula parameters are determined sequentially to reproduce the prescribed correlations in the original variable space using the correlation-matching procedure described above. The resulting moment estimations are presented in Table~\ref{table:ch4-horizontal-truss-gauss-results}, and the corresponding relative errors are shown in Figure~\ref{fig:ch4-horizontal-truss-copula}(a).

The results under the Gaussian copula assumption show that the QPEM provides the most accurate overall estimation among the considered methods. In particular, the relative errors in skewness and kurtosis are approximately 9.1\% and 16.0\%, respectively, while substantially larger deviations are observed for most of the competing methods, especially for the higher-order moments. Nevertheless, the remaining discrepancies in the QPEM estimates indicate the influence of the assumed dependence model. Although the Gaussian copula is a natural baseline when the available dependence information is limited to the correlation matrix, it represents only one of the possible joint dependence structures consistent with the prescribed marginals and correlations. The effect of adopting an alternative non-Gaussian dependence structure is therefore investigated next.

\begin{table}[!t]
	\caption{\label{table:ch4-horizontal-truss-gumbel-results}Moment estimations for the maximum vertical displacement of the horizontal truss structure modeled by Gumbel copula ($n=10$)}
	\begin{tabular*}{\textwidth}[t]{p{0.15\textwidth} p{0.2\textwidth} p{0.14\textwidth} p{0.12\textwidth} p{0.12\textwidth} p{0.12\textwidth}}
		\toprule
		\textbf{Method}   &\centering \textbf{No. of Points} & \textbf{Mean [mm]}& \textbf{STD [mm]} & \textbf{Skewness} &\textbf{ Kurtosis}\\\midrule
		MC     &\centering  $10^6$   & 8.1597   & 2.3359		& 1.1364 & 5.5654 \\ 
		LHS    &\centering  201      & 8.1251     & 2.3959       & 1.3107        & 4.3674\\
		QMC    &\centering  201       & 8.1841      & 2.2069        & 0.9186    & 4.7654 \\
		SGH3   &\centering  221        & 8.1128      & 2.4346        & 0.7051        & 1.5749 \\
		HPEM   &\centering  21        &8.1420      & 2.1395        & -0.3069        & 1.0019 \\
		QPEM ($r=3$)&\centering  201  & 8.1375      & 2.3710        & 1.0863        & 4.9735 \\\bottomrule
	\end{tabular*}
\end{table}

Alternatively, the Gumbel copula family is investigated as a substitute for the Gaussian assumption. For this purpose, a D-vine with the ordering $P_1,\ldots,P_6$ is considered, and all pair-copulas are modeled using the Gumbel family. The corresponding parameters are determined sequentially to reproduce the prescribed correlations among the load variables, while the remaining input variables are maintained as statistically independent. The resulting moment estimations are listed in Table~\ref{table:ch4-horizontal-truss-gumbel-results}, and the relative errors are shown in Figure~\ref{fig:ch4-horizontal-truss-copula}(b).

The computational results indicate that the assumed dependence model has a noticeable influence on the estimated response moments, even though both models are constructed using the same prescribed marginal distributions and correlation matrix. Compared with the Gaussian copula assumption, the Gumbel model improves the QPEM estimates across all four moments, with the most notable improvements observed for the higher-order moments. Specifically, the relative error in skewness is reduced from approximately 9.1\% to 4.4\%, while the kurtosis error decreases from approximately 16.0\% to 10.6\%. The errors in the mean and standard deviation are also slightly reduced. These results demonstrate that different dependence structures consistent with the same prescribed second-order information can lead to different response moment estimations, particularly for higher-order moments.

Overall, the QPEM consistently provides the most accurate moment estimations among the considered methods under both the Gaussian and Gumbel dependence assumptions, achieving the lowest relative errors across all four response moments.

\subsubsection{Example 6: Rock tunnel excavation}
\begin{figure}[!t]
    \centering
    \includegraphics[width=0.6\textwidth]{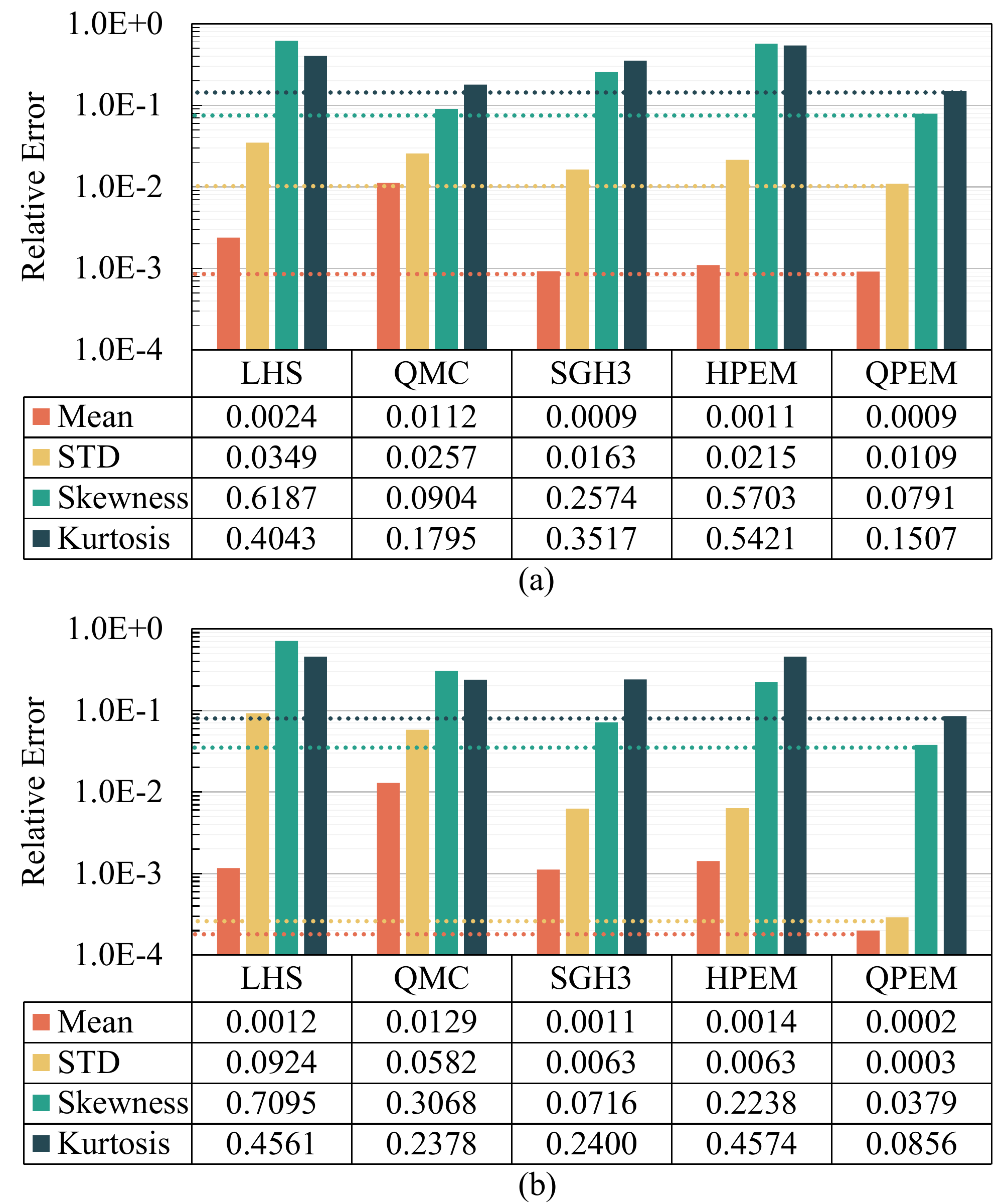}
    \caption{Relative errors for the rock tunnel example under different dependence models: (a) Gaussian copula; and (b) prescribed vine structure and pair-copula families.}
\label{fig:ch4-rock-tunnel-copula}
\end{figure}

The rock tunnel excavation problem introduced previously in Example 2 is revisited to further investigate the effect of dependence-model assumptions under limited probabilistic information. As in the previous example, the available quantitative information is assumed to be limited to the marginal distributions and the correlation matrix. A Gaussian copula, corresponding to the conventional Nataf transformation, is first adopted as a baseline dependence model. The Gaussian-space correlation parameters are determined to reproduce the prescribed correlations in the original variable space. The resulting moment estimations are presented in Table~\ref{table:ch4-rock-tunnel-gauss-results}, and the corresponding relative errors are shown in Figure~\ref{fig:ch4-rock-tunnel-copula}(a).

However, the prescribed marginals and correlations alone do not uniquely determine the joint dependence structure. In many engineering applications, additional qualitative information regarding the dependence structure may also be available from physical understanding, prior knowledge, or engineering judgment. To represent such a scenario, the vine structure and pair-copula families are additionally assumed to be known, while their corresponding parameters remain unknown. The pair-copula parameters are determined using the sequential procedure described above to reproduce the same prescribed correlation matrix. The resulting moment estimations are presented in Table~\ref{table:ch4-rock-tunnel-exact-results}, with the corresponding relative errors shown in Figure~\ref{fig:ch4-rock-tunnel-copula}(b). 

The results show that incorporating the additional qualitative dependence information leads to a noticeable improvement in the QPEM estimates. In particular, the relative error in skewness is reduced from approximately 7.9\% under the Gaussian copula assumption to 3.8\% when the prescribed vine structure and pair-copula families are considered. The kurtosis error is also reduced from approximately 15.1\% to 8.6\%, while further improvements are observed for the mean and standard deviation. Although the estimation of kurtosis remains more challenging, the QPEM provides the lowest relative errors across all four response moments under the informed vine-copula model. These results indicate that, when credible qualitative information regarding the dependence structure is available, incorporating such information can improve probabilistic moment predictions beyond those obtained from the conventional Gaussian dependence assumption.


\begin{table}[!t]
	\caption{\label{table:ch4-rock-tunnel-gauss-results}Moment estimations for the rock tunnel using the Gaussian copula ($n=5$)}

	\begin{tabular*}{\textwidth}[t]{p{0.17\textwidth} p{0.2\textwidth} p{0.12\textwidth} p{0.12\textwidth} p{0.12\textwidth} p{0.12\textwidth}}
		\toprule
		\textbf{Method}   &\centering \textbf{No. of Points} & \textbf{Mean} & \textbf{STD} & \textbf{Skewness} & \textbf{Kurtosis}\\\midrule
		MC     &\centering  $10^6$   & 0.0059  & 0.0009	& -1.0389 & 4.8973 \\ 
		LHS    &\centering  51      & 0.0058     & 0.0009       & -0.3961        & 2.9173 \\
		QMC    &\centering  51       & 0.0058      & 0.0010        & -0.9450        & 4.0182 \\
		SGH3   &\centering  61        & 0.0059       & 0.0009        & -0.7715       & 3.1748 \\
		HPEM   &\centering  11        &0.0059      & 0.0009        & -0.4465        & 2.2427 \\
		QPEM ($r=3$)&\centering  51  & 0.0059      & 0.0009       & -0.9567        & 4.1590 \\\bottomrule
	\end{tabular*}
\end{table}

\begin{table}[!t]
	\caption{\label{table:ch4-rock-tunnel-exact-results}Moment estimations for the rock tunnel using the prescribed vine structure and pair-copula families ($n=5$)}
	\begin{tabular*}{\textwidth}[t]{p{0.17\textwidth} p{0.2\textwidth} p{0.12\textwidth} p{0.12\textwidth} p{0.12\textwidth} p{0.12\textwidth}}
		\toprule
		\textbf{Method}   &\centering \textbf{No. of Points} & \textbf{Mean}& \textbf{STD} & \textbf{Skewness} & \textbf{Kurtosis}\\\midrule
		MC     &\centering  $10^6$   & 0.0059  & 0.0009	& -1.0389 & 4.8973 \\ 
		LHS    &\centering  51      & 0.0058     & 0.0009       & -0.3018        & 2.6638 \\
		QMC    &\centering  51       & 0.0058      & 0.0010        & -1.3576        & 6.0616 \\
		SGH3   &\centering  61        & 0.0059       & 0.0009        & -0.9645       & 3.7220 \\
		HPEM   &\centering  11        &0.0059      & 0.0009        & -0.8064        & 2.6573 \\
		QPEM ($r=3$)&\centering  51  & 0.0059      & 0.0009       & -0.9996        & 4.4780 \\\bottomrule
	\end{tabular*}
\end{table}

\bmsection{Conclusion}
\label{sec:conclusion}

This study presents a unified probabilistic framework that integrates the computational efficiency of the Quadratic Point Estimate Method (QPEM) with the modeling flexibility of vine copulas for uncertainty propagation involving dependent non-Gaussian inputs. The proposed framework employs the inverse Rosenblatt transformation to map the deterministic QPEM sigma points from the independent standard Gaussian space to the original dependent input space. Two representative levels of probabilistic information are considered. In Case (1), multivariate observations are available and the vine structure, pair-copula families, and associated parameters are inferred from data using a sequential heuristic procedure. In Case (2), the available information is restricted to the marginal distributions and correlation matrix, and alternative vine-copula models are constructed by calibrating their parameters to reproduce the prescribed correlations. Importantly, the latter formulation does not assume that the limited marginal and correlation information uniquely determines the joint distribution, but instead provides a systematic means of investigating the sensitivity of probabilistic moment estimates to the assumed dependence model.
The numerical studies demonstrate that the proposed copula-based QPEM provides consistently accurate estimates of the first four response moments at a relatively low computational cost. For Case (1), separate MC-GT and MC-Fitted references have been employed to distinguish the discrepancy associated with finite-sample input-model inference from the numerical error associated with moment propagation. Across the rock slope, tunnel excavation, horizontal truss, and space truss examples, the QPEM generally provided the most balanced and accurate estimates among the considered methods, particularly for skewness and kurtosis. In the 10-dimensional horizontal truss example, the relative errors with respect to MC-Fitted were approximately 0.2\%, 1.2\%, 4.7\%, and 1.9\% for the mean, standard deviation, skewness, and kurtosis, respectively, using 201 model evaluations. In the 15-dimensional space truss example, the corresponding errors were approximately 0.02\%, 0.07\%, 5.1\%, and 3.6\% using 451 model evaluations. These results indicate that the QPEM retains favorable higher-order moment accuracy as the problem dimension increases, while the distinction between MC-GT and MC-Fitted also highlights the importance of separating input-distribution inference error from numerical moment-propagation error.
The Case (2) studies further demonstrate that the assumed dependence structure can influence response moment estimates even when the same marginal distributions and correlation matrix are prescribed. For the horizontal truss example, replacing the Gaussian-copula baseline with an alternative Gumbel vine reduced the QPEM relative errors in skewness and kurtosis from approximately 9.1\% and 16.0\% to 4.4\% and 10.6\%, respectively. For the rock tunnel example, when additional qualitative information regarding the vine structure and pair-copula families was assumed available, the corresponding QPEM errors in skewness and kurtosis decreased from approximately 7.9\% and 15.1\% under the Gaussian assumption to 3.8\% and 8.6\%. These results do not imply that a particular non-Gaussian dependence model is universally preferable; rather, they demonstrate that different joint distributions consistent with the same limited second-order information can produce different probabilistic predictions, and that credible additional dependence information can improve these predictions when available. Overall, the proposed copula-based QPEM provides an efficient and flexible framework for probabilistic moment estimation under both data-informed and limited-information settings, while explicitly accounting for the role of dependence-model assumptions in uncertainty quantification.

\appendix
\bmsection{Maximum Spanning Tree by Prim's algorithm}
\label{sec:app-mst}
\setcounter{figure}{0}
\setcounter{table}{0}
This section details the implementation of Prim's algorithm (\citealt{cormen2022}), utilized to identify the Maximum Spanning Tree (MST) for the R-vine structure. While Prim's algorithm is conventionally described for finding a \textit{minimum} spanning tree, it operates analogously for maximization problems by simply prioritizing the largest edge weights instead of the smallest.

\begin{figure}[t]
    \centering
    \includegraphics[width=\textwidth]{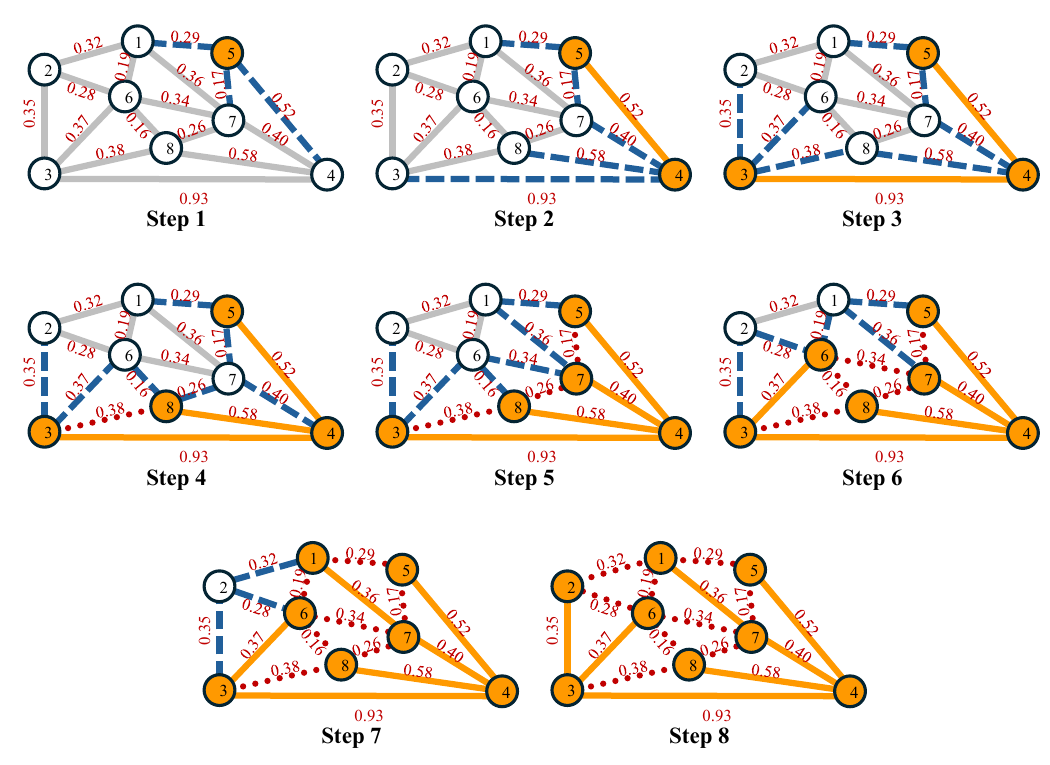}
\caption{An example run through Prim’s Algorithm.}
\label{fig:mst}
\end{figure}

The objective is to find a tree $T=(V, E')$, where $E' \subseteq E$ denotes the set of selected edges, that spans the set of variables $V=\{1, \dots, n\}$ such that the sum of the edge weights (e.g., absolute Kendall's $\tau$) is maximized. The procedure is a greedy algorithm that grows the tree from an arbitrary root node, adding one vertex at a time. To illustrate the algorithm, consider an example graph with 8 nodes ($n=8$). The corresponding evolution of the algorithm is visually depicted in Figure~\ref{fig:mst}, and the general pseudocode is presented in Table~\ref{table:mst_algorithm2}.

\begin{table}[b!]
	\caption{Maximum Spanning Tree (Prim's algorithm)}\label{table:mst_algorithm2}
	\begin{tabular}[t]{p{0.97\textwidth}}
        \hline
	\textbf{Input:} A connected graph $G = (V, E)$ with edge weights $w(e)$.\\
        \textbf{Output:} A set of edges $E'$ forming a Maximum Spanning Tree.\\ [-20pt]
        \begin{itemize}
            \item[\textbf{1}:] Initialize $V_{in} = \{v_{start}\}$, where $v_{start}$ is an arbitrary starting node.
            \item[\textbf{2}:] Initialize $V_{out} = V \setminus \{v_{start}\}$.
            \item[\textbf{3}:] Initialize $E' = \emptyset$.
            \item[\textbf{4}:] \textbf{while} $V_{out} \neq \emptyset$ \textbf{do}
            \item[\textbf{5}:] \qquad Find the edge $e^* = (u, v)$ with maximal weight such that $u \in V_{in}$ and $v \in V_{out}$.
            \item[\textbf{6}:] \qquad $V_{in} \leftarrow V_{in} \cup \{v\}$
            \item[\textbf{7}:] \qquad $V_{out} \leftarrow V_{out} \setminus \{v\}$
            \item[\textbf{8}:] \qquad $E' \leftarrow E' \cup \{e^*\}$
            \item[\textbf{9}:] \textbf{end while}
            \item[\textbf{10}:] \textbf{return} $E'$
        \end{itemize}\\ 
        \hline
	\end{tabular}
\end{table}

The step-by-step execution of the algorithm proceeds as follows:

\begin{itemize}
    \item \textbf{Initialization:} The algorithm begins by selecting an arbitrary starting node. In this example, node 5 is chosen as the root. The set of visited nodes is initialized as $V_{in}=\{5\}$, and the set of unvisited nodes is $V_{out}=\{1, \dots, 8\} \backslash \{5\}$.

    \item \textbf{Step 1:} We examine all edges connecting the node in $V_{in}$ to any node in $V_{out}$ (blue dashed lines). Among all candidate edges, the edge with the maximum absolute Kendall's $\tau$ is selected. As shown in Figure~\ref{fig:mst}, the edge $(5, 4)$ is chosen. Node 4 is added to the tree, updating $V_{in}=\{4, 5\}$.
    
    \item \textbf{Steps 2-3:} In each subsequent step, the algorithm considers all possible edges (blue dashed lines) connecting the current set of visited nodes $V_{in}$ to the remaining unvisited nodes $V_{out}$. The edge with the highest weight is consistently selected and added to the spanning tree. This process ensures that the strongest dependencies are prioritized. The specific edges selected at each step and their corresponding weights are listed in Table~\ref{table:mst_algorithm1}.

    \item \textbf{Step 4:} Since both nodes 3 and 8 are already in $V_{in}$, adding this edge would create a \textit{cycle} (closed loop), contradicting the definition of a tree. Consequently, this edge (3, 8) is identified as invalid and is excluded from the spanning tree, regardless of its edge weight.

    \item \textbf{Steps 5-7:} The procedure repeats. Valid edges with the highest weights are sequentially added to connect the remaining unvisited nodes.
    
    \item \textbf{Step 8:} At the end of Step 7, the edge $(3, 2)$ is added (see Table~\ref{table:mst_algorithm1}). At this point, all 8 nodes are included in $V_{in}$ (i.e., $V_{out} = \emptyset$). Since the spanning tree for $n=8$ nodes must contain exactly $n-1=7$ edges, the algorithm terminates successfully.
\end{itemize}

The resulting set of edges constitutes the MST, which serves as the first tree ($T_1$) in the R-vine structure. This structure captures the most significant pairwise dependencies among the 8 variables.
\begin{table}[t!]
\centering
\caption{Step-by-step execution of Prim's algorithm for MST construction ($n=8$)}
\label{table:mst_algorithm1}
\begin{tabular}{p{0.05\textwidth} p{0.23\textwidth} p{0.27\textwidth} p{0.2\textwidth} p{0.13\textwidth}}
\hline
\centering \textbf{Step} &  \centering \textbf{Visited Nodes $V_{in}$} &  \centering \textbf{Unvisited Nodes $V_{out}$} &\centering \textbf{Selected Edge $e^*$} & \quad \; \textbf{Weight} \\ 
\hline
\centering 1 & \centering $\{5\}$ & \centering $\{1,2,3, 4, 6, 7, 8\}$ & \centering $(5,4)$ & \qquad \; 0.52 \\
\centering 2 & \centering $\{4, 5\}$ & \centering $\{1,2,3, 6, 7, 8\}$ &  \centering $(4,3)$ & \qquad \; 0.93 \\
\centering 3 & \centering $\{3, 4, 5\}$ & \centering $\{1,2, 6, 7, 8\}$ &  \centering $(4,8)$ & \qquad \; 0.58 \\
\centering 4 & \centering $\{3, 4, 5, 8\}$ & \centering $\{1,2,6, 7\}$ &  \centering $(4,7)$ & \qquad \; 0.40 \\
\centering 5 & \centering $\{3, 4, 5, 7, 8\}$ & \centering $\{1,2,6\}$ &  \centering $(3,6)$ & \qquad \; 0.37 \\
\centering 6 & \centering $\{3, 4, 5, 6, 7, 8\}$ & \centering $\{1,2\}$ &  \centering $(7, 1)$ & \qquad \; 0.36 \\
\centering 7 & \centering $\{1,3, 4, 5, 6, 7, 8\}$ & \centering $\{2\}$ &  \centering $(3,2)$ & \qquad \; 0.35 \\
\centering 8 & \centering $\{1,2,3, 4, 5, 6, 7, 8\}$ & \centering $\emptyset$ & \centering -- & \qquad \;\;\; -- \\
\hline
\end{tabular}
\end{table}
Notice that we did not use the actual values of the edges here, instead we only used their rank. Therefore, the algorithm leads to the same results if we transform the edge values by a monotone increasing function. Hence, in our area of application, where we want to find a tree with maximal values of absolute Kendall’s $\tau$, we would get the same tree even if we consider other weights, like squared $\tau$ or another monotone increasing function.

\bmsection{Parameter Estimation and Model Selection}
\label{app:estimation}

While the Maximum Likelihood Estimator (MLE) theoretically provides the most efficient estimates for vine copula parameters, its practical application in high-dimensional engineering problems may entail a prohibitive computational burden. The simultaneous optimization of parameters across all tree levels often leads to numerical instability and intractability. To circumvent these challenges, this study adopts the \textit{Stepwise Semiparametric (SSP)} estimator (\citealt{haff2013}). This method effectively decouples the estimation of marginal distributions from the dependence structure and further decomposes the high-dimensional copula estimation into a sequence of low-dimensional optimization problems.

The SSP estimation procedure is rigorously defined in the following two stages:

\subsection*{\textbf{Stage 1: Transformation to Pseudo-observations via Empirical Marginals}}
\noindent The Empirical Cumulative Distribution Function (ECDF) is employed to handle the marginals non-parametrically. This step transforms the raw data into uniform variates.

Let $\mathbf{X} = (\mathbf{x}_1, \dots, \mathbf{x}_N)^T$ denote the data matrix containing $N$ independent realizations. The observations are transformed into \textit{pseudo-observations} $\mathbf{U}$ based on their empirical ranks:
\begin{equation}
    \hat{u}_{\ell i} = \hat{F}_i(x_{\ell i}) = \frac{1}{N+1} \sum_{k=1}^N \mathbb{I}(x_{k i} \le x_{\ell i}) = \frac{R_{\ell i}}{N+1}
\end{equation}
where $\mathbb{I}(\cdot)$ is the indicator function and $R_{\ell i}$ denotes the rank of $x_{\ell i}$ among the observations $\{x_{1i}, \dots, x_{Ni}\}$. The divisor $N+1$ ensures that the transformed values strictly remain within the open interval $(0,1)$. These normalized ranks serve as the input for the subsequent copula estimation (\citealt{genest1995semiparametric}).

\subsection*{\textbf{Stage 2: Sequential (Stepwise) Estimation of Copula Parameters}}
\noindent Given the pseudo-observations, the copula parameters are estimated sequentially, tree by tree. This approach relies on the hierarchical structure of the R-vine, where the conditional distributions used in tree $T_j$ are derived from the pair-copulas estimated in trees $T_1, \dots, T_{j-1}$.

Let $\boldsymbol{\theta}^k$ represent the parameter vector for all pair-copulas in the $k$-th tree. The SSP estimator obtains $\hat{\boldsymbol{\theta}}^k$ by maximizing the pseudo-log-likelihood function specific to that level:

\begin{itemize}
    \item \textbf{Level 1 ($k=1$):} The parameters for the first tree are estimated directly from the pseudo-observations:
    \begin{equation}
        \hat{\boldsymbol{\theta}}^1 = \operatorname*{argmax}_{\boldsymbol{\theta}^1} \sum_{\ell=1}^N \sum_{e \in E_1} \ln c_{e}\left(\hat{u}_{\ell, i_e}, \hat{u}_{\ell, j_e}; \theta_e\right)
    \end{equation}
    where $E_1$ is the edge set of $T_1$, and $(i_e, j_e)$ are the indices of the variables connected by edge $e$.

    \item \textbf{Level $k$ ($k \ge 2$):} The estimation for subsequent trees requires transformed observations. Let $v_{e,1}^{(\ell)}$ and $v_{e,2}^{(\ell)}$ denote the conditional pseudo-observations for edge $e \in E_k$, computed recursively using the $h$-functions (Eq.~\ref{eq:h-function}) with the estimated parameters $\hat{\boldsymbol{\theta}}^{1:k-1}$ from previous levels. The parameters for tree $T_k$ are then obtained by:
    \begin{equation}
        \hat{\boldsymbol{\theta}}^k = \operatorname*{argmax}_{\boldsymbol{\theta}^k} \sum_{\ell=1}^N \sum_{e \in E_k} \ln c_{e}\left(v_{e,1}^{(\ell)}, v_{e,2}^{(\ell)}; \theta_e\right)
    \end{equation}
\end{itemize}
Although the SSP estimator is asymptotically less efficient than the global MLE due to the stepwise accumulation of estimation errors, \cite{haff2013} demonstrated that it remains consistent and asymptotically normal, providing a rigorous yet computationally feasible solution for high-dimensional inference.

\subsection{Model Selection Criterion}
Within the SSP framework, the specific parametric family for each pair-copula is selected to best represent the local dependence structure. For every edge $e$ in the vine, a set of candidate families $\mathcal{M}$ (e.g., Gaussian, Student-$t$, Clayton, Gumbel, Frank, etc.) is evaluated. 

While the Bayesian Information Criterion (BIC) constitutes a valid alternative, imposing a stricter penalty for model complexity, this study employs the Akaike Information Criterion (AIC). The AIC is favored for its asymptotic efficiency in minimizing the information loss relative to the true distribution (\citealt{burnham2002model}). The optimal family $m^* \in \mathcal{M}$ is identified by minimizing the AIC defined as:
\begin{equation}
    \text{AIC}_{e} = -2 \sum_{\ell=1}^N \ln c_{e}\left(u_{1}^{(\ell)}, u_{2}^{(\ell)}; \hat{\theta}_e\right) + 2 p_m
\end{equation}
where $p_m$ denotes the number of parameters of the copula family $m$. This criterion ensures a balance between goodness-of-fit and parsimony in the final vine copula model.

\bmsection{Definitions of Pair-Copula Families}
\label{sec:app-copula-def}

This Appendix provides the cumulative distribution functions (CDFs) of the bivariate copula families considered in this study, followed by the formulation for rotated copulas used to model negative dependencies.

\subsection{Formulas of Copula Families}
Let $u$ and $v$ denote the uniform marginals in $[0,1]$. The formulas for the copula families $C(u,v)$ and their parameter ranges are listed below (\citealt{nelsen2007copula,joe2014dependence, torre2019general}).

\subsubsection*{Gaussian Copula}
\begin{equation}
    C(u, v; \theta) = \Phi_2\left(\Phi^{-1}(u), \Phi^{-1}(v); \theta\right)
\end{equation}
where $\Phi_2(\cdot, \cdot; \theta)$ is the bivariate standard normal CDF with dependence parameter $\theta \in (-1, 1)$, and $\Phi^{-1}$ is the inverse univariate standard normal CDF.

\subsubsection*{Student's t-Copula}
\begin{equation}
    C(u, v; \theta, \nu) = t_{2, \nu}\left(t_{\nu}^{-1}(u), t_{\nu}^{-1}(v); \theta, \nu\right)
\end{equation}
where $t_{2, \nu}$ is the bivariate Student's t-distribution CDF with dependence parameter $\theta \in (-1, 1)$ and degrees of freedom $\nu > 1$, and $t_{\nu}^{-1}$ is the inverse univariate Student's t-distribution CDF.

\subsubsection*{Clayton Copula}
\begin{equation}
    C(u, v; \theta) = \left( \max\left\{ u^{-\theta} + v^{-\theta} - 1, 0 \right\} \right)^{-1/\theta}, \quad \theta > 0
\end{equation}

\subsubsection*{Gumbel Copula}
\begin{equation}
    C(u, v; \theta) = \exp\left( -\left[ (-\ln u)^{\theta} + (-\ln v)^{\theta} \right]^{1/\theta} \right), \quad \theta \ge 1
\end{equation}

\subsubsection*{Frank Copula}
\begin{equation}
    C(u, v; \theta) = -\frac{1}{\theta} \ln\left( 1 + \frac{(e^{-\theta u} - 1)(e^{-\theta v} - 1)}{e^{-\theta} - 1} \right), \quad \theta \in \mathbb{R} \backslash \{0\}
\end{equation}

\subsubsection*{Joe Copula}
\begin{equation}
    C(u, v; \theta) = 1 - \left( (1-u)^{\theta} + (1-v)^{\theta} - (1-u)^{\theta}(1-v)^{\theta} \right)^{1/\theta}, \quad \theta \ge 1
\end{equation}

\subsubsection*{Ali-Mikhail-Haq Copula}
\begin{equation}
    C(u, v; \theta) = \frac{uv}{1 - \theta(1-u)(1-v)}, \quad \theta \in [-1, 1]
\end{equation}

\subsubsection*{Farlie-Gumbel-Morgenstern Copula}
\begin{equation}
    C(u, v; \theta) = uv\left( 1 + \theta(1-u)(1-v) \right), \quad \theta \in [-1, 1] \backslash \{0\}
\end{equation}

\subsubsection*{Plackett Copula}
\begin{equation}
    C(u, v; \theta) = \frac{1}{2(\theta-1)} \left( 1+(\theta-1)(u+v) - \sqrt{\left[1+(\theta-1)(u+v)\right]^2 - 4\theta(\theta-1)uv} \right), \quad \theta > 0, \theta \ne 1
\end{equation}

\subsection{Rotated Copulas}
Standard Archimedean copulas (e.g., Clayton, Gumbel, Joe) typically capture only positive dependence. To model negative dependence or different tail behaviors, these copulas can be rotated. The rotated copulas are obtained by transforming the arguments $u$ and $v$ as follows:

\begin{itemize}
    \item \textbf{Rotation by 90$^{\circ}$}:
    \begin{equation}
        C^{90^{\circ}}(u, v) = v - C(1-u, v)
    \end{equation}
    This rotation captures negative dependence. For example, a 90$^{\circ}$-rotated Clayton copula exhibits dependence in the lower-right tail.
    
    \item \textbf{Rotation by 180$^{\circ}$ (Survival Copula)}:
    \begin{equation}
        C^{180^{\circ}}(u, v) = u + v - 1 + C(1-u, 1-v)
    \end{equation}
    This is often referred to as the survival copula. It reflects the dependence structure of the standard copula rotated by 180 degrees (e.g., upper-tail dependence becomes lower-tail dependence).
    
    \item \textbf{Rotation by 270$^{\circ}$}:
    \begin{equation}
        C^{270^{\circ}}(u, v) = u - C(u, 1-v)
    \end{equation}
    Similar to the 90$^{\circ}$ rotation, this configuration captures negative dependence, typically in the upper-left tail.
\end{itemize}

For symmetric copulas such as Frank, Gaussian, and Student's t, rotations are related to sign changes in their parameters (e.g., $C_{\text{Frank}, 90^{\circ}}(u,v;\theta) = C_{\text{Frank}}(u,v;-\theta)$). However, for asymmetric copulas (Clayton, Gumbel, Joe), these rotational transformations provide distinct dependence structures crucial for flexible modeling.

\bmsection{Target Correlation Matrices for Case (2)}
\label{app:target_corr}
For completeness, the target correlation matrices employed in the Case (2) numerical examples are provided below. These matrices represent the prescribed second-order dependence information used for constructing the corresponding input probabilistic models. Since closed-form Pearson correlations are generally unavailable for the prescribed non-Gaussian vine-copula models, the target correlations are evaluated numerically from the corresponding ground-truth input distributions. These calculations involve only input-distribution sampling and do not require evaluations of the computational response models.

\subsection{Example (5): Maximum vertical displacement of horizontal truss}
\begin{equation}
    \mathbf{R}_{\mathrm{target}} =
    \begin{bmatrix}
        1.00 & 0.00 & 0.00 & 0.00 & 0.00 & 0.00 & 0.00 & 0.00 & 0.00 & 0.00 \\
        0.00 & 1.00 & 0.00 & 0.00 & 0.00 & 0.00 & 0.00 & 0.00 & 0.00 & 0.00 \\
        0.00 & 0.00 & 1.00 & 0.00 & 0.00 & 0.00 & 0.00 & 0.00 & 0.00 & 0.00 \\
        0.00 & 0.00 & 0.00 & 1.00 & 0.00 & 0.00 & 0.00 & 0.00 & 0.00 & 0.00 \\
        0.00 & 0.00 & 0.00 & 0.00 & 1.00 & 0.75 & 0.84 & 0.80 & 0.76 & 0.77 \\
        0.00 & 0.00 & 0.00 & 0.00 & 0.75 & 1.00 & 0.75 & 0.84 & 0.80 & 0.76 \\
        0.00 & 0.00 & 0.00 & 0.00 & 0.84 & 0.75 & 1.00 & 0.75 & 0.84 & 0.80 \\
        0.00 & 0.00 & 0.00 & 0.00 & 0.80 & 0.84 & 0.75 & 1.00 & 0.75 & 0.84 \\
        0.00 & 0.00 & 0.00 & 0.00 & 0.76 & 0.80 & 0.84 & 0.75 & 1.00 & 0.75 \\
        0.00 & 0.00 & 0.00 & 0.00 & 0.77 & 0.76 & 0.80 & 0.84 & 0.75 & 1.00
    \end{bmatrix}
\label{eq:target_corr_horizontal_truss}
\end{equation}

\subsection{Example (6): Rock tunnel excavation}
\begin{equation}
    \mathbf{R}_{\mathrm{target}} =
    \begin{bmatrix}
        1.00 & 0.51 & 0.52 & -0.24 & 0.33 \\
        0.51 & 1.00 & 0.13 & -0.13 & 0.21 \\
        0.52 & 0.13 & 1.00 & -0.13 & 0.21 \\
        -0.24 & -0.13 & -0.13 & 1.00 & -0.08 \\
        0.33 & 0.21 & 0.21 & -0.08 & 1.00
    \end{bmatrix}
\label{eq:target_corr_tunnel}
\end{equation}

\bmsectionstar{Data Availability Statement}

Some or all data, models, or code that support the findings of this study are available from the corresponding author upon reasonable request.

\stepcounter{pdfbookmarksec}
\pdfbookmark[1]{References}{bmsec.\arabic{pdfbookmarksec}}
\bibliography{reference_copula}%

\end{document}